%% file: paper.tex
\documentclass[a4paper]{article}
\usepackage[inner=30mm,outer=30mm,textheight=225mm]{geometry}

\newcommand{\onlyau}[1]{}

\usepackage{amsmath}
\usepackage{amssymb}
\usepackage{enumerate}
\usepackage{graphicx}
\usepackage{psfrag}
\usepackage{supertabular}
\usepackage{theorem}

\newtheorem{theorem}{Theorem}[section]
\newtheorem{corollary}[theorem]{Corollary}
\newtheorem{defn}[theorem]{Definition}
\newtheorem{lemma}[theorem]{Lemma}
\newtheorem{algorithm}[theorem]{Algorithm}

\newcommand{\prooflabel}{Proof}
\newcommand{\qed}{\hfill\rule[-0.5mm]{1.5mm}{3.0mm}}
\newtheorem{proofthm}{\prooflabel}
\newenvironment{proof}{\begin{proofthm} \em}{\qed \end{proofthm}}

\newcommand{\annhalfref}{\dot{A}}
\newcommand{\annref}{\bar{A}}
\newcommand{\annrefref}{\bar{\bar{A}}}
\newcommand{\dischalfref}{\dot{D}}
\newcommand{\discref}{\bar{D}}
\newcommand{\edge}[1]{\overline{#1}}
\newcommand{\dedge}[1]{\overrightarrow{#1}}
\newcommand{\kb}{K^2}
\newcommand{\mb}{M}
\newcommand{\mobius}{M\"{o}bius}
\newcommand{\mobref}{\bar{M}}
\newcommand{\polyedge}[1]{$\edge{\mathit{#1}}$}
\newcommand{\polydedge}[1]{$\dedge{\mathit{#1}}$}
\newcommand{\polytri}[1]{$\Delta \mathit{#1}$}
\newcommand{\ppirr}{{$\mathbb{P}^2$-ir\-re\-du\-ci\-ble}}
\newcommand{\ppirrty}{{$\mathbb{P}^2$-ir\-re\-du\-ci\-bi\-li\-ty}}
\newcommand{\regina}{{\em Regina}}
\newcommand{\rpp}{\mathbb{R}P^2}
\newcommand{\rps}{\mathbb{R}P^3}
\newcommand{\sfs}[2]{\mathrm{SFS}\left(#1: #2\right)}
\newcommand{\sfslong}{Seifert fibred space}

\newcommand{\snappea}{{\em SnapPea}}
\newcommand{\specialnine}{\mathrm{SFS}\left(\smash{\annhalfref}/o_2\right) /
  {\homtwotable{0}{1}{1}{0}}}
\newcommand{\specialten}{\sfs{\smash{\dischalfref}}{(2,1)} \allowbreak
  \bigcup \allowbreak \mathit{Gieseking}}

\newcommand{\sss}{S^3}
\newcommand{\torus}{T^2}
\newcommand{\Z}{\mathbb{Z}}

\newcommand{\homtwotable}[4]{
    \mbox{\tiny \renewcommand{\arraystretch}{1}
        $\! \left[ \begin{array}{@{\ }r@{\ }r@{\ }} #1 & #2 \\ #3 & #4
        \end{array} \right]$
    }
}

\newcommand{\twolines}[3]{
    \begin{tabular}{#1} #2 \\ #3 \end{tabular}
}

\title{Enumeration of non-orientable 3-manifolds \\
    using face pairing graphs and union-find}
\author{Benjamin A.~Burton}

\date{April 27, 2006}

\begin{document}

\maketitle

\abstract{Drawing together techniques from combinatorics and
    computer science, we improve the census algorithm for enumerating
    closed minimal {\ppirr} 3-manifold triangulations.  In particular,
    new constraints are proven for face pairing graphs, and
    pruning techniques are improved using a modification of the
    union-find algorithm.  Using these results we catalogue all 136 closed
    non-orientable {\ppirr} 3-manifolds that can be formed from at most
    ten tetrahedra.}

\input{intro.tex}

\input{algm.tex}

\input{prelim.tex}

\input{graphs.tex}

\input{ufind.tex}

\input{census.tex}

\section*{Acknowledgements}

The author is grateful to the Victorian Partnership for Advanced
Computing for supporting this work under e-Research grant EPANRM155.2005,
and also for the use of their exceptional computational resources.

\bibliographystyle{amsplain}
\bibliography{pure-20060425}

\vspace{1cm}
\noindent
Benjamin A.~Burton \\
Department of Mathematics, SMGS, RMIT University \\
GPO Box 2476V, Melbourne, VIC 3001, Australia \\
(bab@debian.org)

\appendix
\input{manifolds.tex}

\end{document}

%% file: intro.tex
\section{Introduction}

With recent advances in computing power, topologists have been able to
construct exhaustive tables of small 3-manifold triangulations, much like
knot theorists have constructed exhaustive tables of simple knot projections.
Such tables are valuable sources of data, but they suffer from the fact
that enormous amounts of computer time are required to generate them.

Where knot tables are often limited by bounding the number of
crossings in a knot projection, tables of 3-manifolds generally limit
the number of tetrahedra used in a 3-manifold triangulation.  A typical
table (or {\em census}) of 3-manifolds lists all 3-manifolds of a
particular type that can be formed from $t$ tetrahedra or fewer.

Beyond their generic role as a rich source of examples, tables of
this form have a number of specific uses.
They still offer the only general means for proving that a triangulation
is minimal (i.e., uses as few tetrahedra as possible), much in the same way
as knot tables
are used to calculate crossing number.  Moreover, a detailed analysis of
these tables can offer insight into the combinatorial structures of
minimal triangulations, as seen for example by the structural observations
of Matveev \cite{matveev6}, Martelli and Petronio \cite{italianfamilies}
and Burton \cite{burton-nor7}.

Unfortunately the scope of such tables is limited by the difficulty of
generating them.  In general, a census of triangulations formed from
$\leq t$ tetrahedra requires computing time at least exponential in $t$.
In the case of closed 3-manifold triangulations, results are only known
for $\leq 11$ tetrahedra in the orientable case and $\leq 8$ tetrahedra
in the non-orientable case.  These results are particularly sparse in
the non-orientable case --- only 18 distinct manifolds are found, all of
which are graph manifolds \cite{burton-nor8}.

Clearly there is more to be learned by extending the existing censuses
to higher numbers of tetrahedra.  Due to the heavy computational
requirements however, this requires significant improvements in the
algorithms used to generate the census data.  Such improvements form the
main subject of this paper.

We restrict our attention here to closed 3-manifold triangulations.
In the orientable case, successive tables have been generated by
Matveev \cite{matveev6}, Ovchinnikov,
Martelli and Petronio \cite{italian9}, Martelli \cite{italian10} and
then Matveev again with 11~tetrahedra \cite{matveev11}.
For non-orientable manifolds, tabulation begins with
Amendola and Martelli \cite{italian-nor6,italian-nor7} and is continued
by Burton \cite{burton-nor7,burton-nor8} up to 8~tetrahedra.

The contributions of this paper are the following:
\begin{itemize}
    \item Several improvements to the algorithm for generating
    census data, some of which increase the speed by orders of
    magnitude;
    \item An extension of the closed non-orientable census from
    8~tetrahedra to 10~tetrahedra;
    \item A verification of previous closed orientable census results
    for up to 10~tetrahedra, and an extension of these results from a census
    of manifolds to a census of all minimal triangulations.
\end{itemize}

The algorithmic improvements are divided into two broad categories.
The first set of results relate to {\em face pairing graphs}, which are
4-valent graphs describing which tetrahedron faces are identified to
which within a triangulation.  The second set is based upon
the {\em union-find} algorithm, which is a well-known method for
finding connected components in a graph.  Here the union-find
algorithm is modified to support backtracking, and to efficiently
monitor properties of vertex and edge links within a triangulation.

All computational work was performed using the topological
software package {\regina} \cite{regina,burton-regina}.  Census
generation forms only a small part of {\regina}, which is a larger
software package for performing a variety of tasks in 3-manifold
topology.
The software is released under the GNU General Public License, and
may be freely downloaded from
{\tt http://\allowbreak regina.\allowbreak sourceforge.\allowbreak net/}.

The bulk of this paper is devoted to improving the census algorithm.
Section~\ref{s-algm} begins with the precise census constraints, and
follows with an overview of how a census algorithm is structured.
In Section~\ref{s-prelim} we present a series of preliminary results,
describing properties of minimal triangulations that will be required
in later sections.  Section~\ref{s-graphs} offers the first
round of algorithmic improvements, based upon the analysis of face pairing
graphs.  A more striking set of improvements is made in
Section~\ref{s-ufind}, in which a modified union-find algorithm is used
to greatly reduce the search space.  Both Sections~\ref{s-graphs}
and~\ref{s-ufind} also include empirical results in which the effectiveness
of these improvements is measured.

The improvements of Sections~\ref{s-graphs} and~\ref{s-ufind}
have led to new closed census results, as outlined
above.  Section~\ref{s-census} summarises these new results,
with a focus on the extension of the closed non-orientable census from
8 to 10 tetrahedra.  A full list of non-orientable census manifolds is
included in the appendix.

%% file: algm.tex
\section{Overview of the Census Algorithm} \label{s-algm}

As is usual for a census of closed 3-manifolds, we restrict our attention to
manifolds with the following properties:
\begin{itemize}
    \item {\em Closed:} The 3-manifold is compact, with no boundary and no
    cusps.
    \item {\em {\ppirr}:} The 3-manifold contains no embedded two-sided
    projective planes, and every embedded 2-sphere bounds a ball.
\end{itemize}

The additional constraint of {\ppirrty} allows us to focus on the most
``fundamental'' manifolds --- the properties of larger manifolds are often
well understood in terms of their {\ppirr} constituents.

Recall that we are not just enumerating 3-manifolds, but also their
triangulations.  We therefore focus only on triangulations with the
following additional property:
\begin{itemize}
    \item {\em Minimal:} The triangulation uses as few tetrahedra as
    possible.  That is, the underlying 3-manifold cannot be triangulated
    using a smaller number of tetrahedra.
\end{itemize}

This minimality constraint is natural for a census, and is used
throughout the literature.
Note that a 3-manifold may have many different minimal triangulations,
though of course all of these triangulations must use the same
number of tetrahedra.

Minimal triangulations are tightly related to the Matveev complexity
of a manifold \cite{matveev-complexity}.  Matveev defines complexity in
terms of special spines, and it has been proven by Matveev in the orientable
case and Martelli and Petronio in the non-orientable case
\cite{italian-decomp} that, with the exceptions of $\sss$, $\rps$ and
$L_{3,1}$, the Matveev complexity of a closed {\ppirr} 3-manifold is
precisely the number of tetrahedra in its minimal triangulation(s).

\subsection{Stages of the Algorithm}

There are two stages involved in constructing a census of 3-manifold
triangulations: the generation of triangulations, and then
the analysis of these triangulations.

\begin{enumerate}
    \item {\em Generation:}
    The generation stage typically involves a long computer search,
    in which tetrahedra are pieced together in all possible ways to
    form 3-manifold triangulations that might satisfy our census constraints.

    The result of this search is a large set of triangulations,
    guaranteed to include all of the triangulations that should be in
    the census.  There are often unwanted triangulations also (for instance,
    triangulations that are non-minimal, or that represent reducible
    manifolds).  This is not a problem; these unwanted triangulations
    will be discarded in the analysis stage.

    The generation of triangulations is entirely automated,
    but it is also extremely time-consuming --- it may take
    seconds or centuries, depending upon the size of the census.

    \item {\em Analysis:}
    Once the generation stage has produced a raw set of triangulations,
    these must be refined into a final census.  This includes verifying
    that each triangulation is minimal and {\ppirr} (and throwing away
    those triangulations that are not).  It also involves grouping
    triangulations into classes that represent the same 3-manifold, and
    identifying these 3-manifolds.

    Analysis is much faster than generation, but it typically requires
    a mixture of automation and human involvement.  Techniques include
    the analysis of invariants and normal surfaces, combinatorial
    analysis of the triangulation structures, and applying elementary
    moves that change triangulations without altering their underlying
    3-manifolds.
\end{enumerate}

The generation stage is the critical bottleneck, due to the vast number
of potential triangulations that can be formed from a small number of
tetrahedra.  Suppose we are searching for triangulations that can be
formed using $n$ tetrahedra.  Even assuming that we know which tetrahedron
faces are to be joined with which, each pair of faces can be identified
according to one of six possible rotations or reflections, giving rise
to $6^{2n}$ possible triangulations in total.  For 10 tetrahedra, this
figure is larger than $10^{15}$.  It is clear then why existing census
data is limited to the small bounds that have been reached to date.

It should be noted that for an orientable census, the figure
$6^{2n}$ becomes closer to $3^n 6^n$.
This is because for a little over half the faces only three of the
six rotations or reflections will preserve orientation.  It is partly
for this reason that the orientable census has consistently been further
advanced in the literature than the non-orientable census.
Nevertheless, for 10 tetrahedra this figure is
still larger than $10^{12}$, a hefty workload indeed.

Fortunately there are ways in which these figures can be reduced.
The usual technique is to use topological arguments to derive
constraints that must be satisfied by minimal {\ppirr} triangulations,
and then to incorporate
these constraints into the generation algorithm so that the search space
can be reduced.  Such techniques have been used throughout the history
of census generation; examples include the edge degree results from the
early cusped hyperbolic census of Hildebrand and Weeks \cite{cuspedcensusold},
and the properties of minimal spines from Matveev's 6-tetrahedron
closed orientable census \cite{matveev6}.

\subsection{Generation of Triangulations} \label{s-algm-gen}

Since generation is the most time-intensive stage of the
census algorithm, it is there that we focus our attention.  Before
describing the generation process in more detail, it is useful to make
the following definition.

\begin{defn}[Face Pairing Graph] \label{d-facegraph}
    Let $T$ be a closed 3-manifold triangulation formed from $n$ tetrahedra.
    The {\em face pairing graph} of $T$ is a 4-valent multigraph
    on $n$ vertices describing
    which faces of which tetrahedra in $T$ are identified.

    More specifically, each vertex of the face pairing graph represents
    a tetrahedron of $T$, and each edge represents a pair of tetrahedron
    faces that are identified.  Note that loops and multiple edges may
    (and frequently do) occur.
\end{defn}

\begin{figure}[htb]
    \centerline{\includegraphics[scale=0.6]{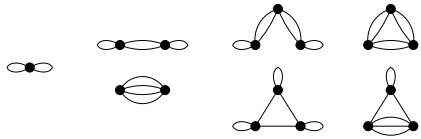}}
    \caption{All possible face pairing graphs for $\leq 3$ tetrahedra}
    \label{fig-graphs3}
\end{figure}

It is clear that every face pairing graph must be 4-valent, since each
tetrahedron has four faces.  Figure~\ref{fig-graphs3} illustrates all
possible face pairing graphs for connected triangulations formed from
$\leq 3$ tetrahedra.

We return now to the generation process.
Algorithm~\ref{a-generation} outlines this process in more detail.
This general scheme has been used by authors throughout the
literature; a similar overview can be found in the early cusped
hyperbolic census of Hildebrand and Weeks \cite{cuspedcensusold}.

\begin{algorithm}[Generation of Triangulations] \label{a-generation}
    Suppose that we wish to generate a set $S$ that includes all closed
    minimal {\ppirr} triangulations formed from precisely $n$ tetrahedra
    (note that $S$ might also include some additional unwanted triangulations).
    The procedure is as follows.
    \begin{enumerate}
        \item {\em Enumerate face pairing graphs:} \label{e-enumgraphs}

        We enumerate (up to isomorphism) all connected 4-valent multigraphs
        on $n$ vertices.  Each of these becomes a candidate face pairing
        graph for an eventual triangulation.

        \item {\em Process face pairing graphs:} \label{e-processgraphs}

        Each candidate face pairing graph generated in
        step~\ref{e-enumgraphs} is processed as follows.  For each edge
        of the graph, we consider all six rotations and reflections
        by which the two corresponding tetrahedron faces might be identified.
        In this way we obtain $6^{2n}$ possible triangulations.

        In practice, constraints that must be satisfied by minimal
        {\ppirr} triangulations are used to prune this search.  As a
        result, far fewer than $6^{2n}$ possible triangulations need to
        be processed.

        Some of these possible triangulations might not be
        triangulations of closed 3-manifolds (e.g., their vertex links
        might not be spheres, or they might have edges that are identified
        with themselves in reverse).
        Each triangulation that is indeed a closed 3-manifold
        triangulation is added to the final set $S$.
    \end{enumerate}
\end{algorithm}

The enumeration of face pairing graphs (step~\ref{e-enumgraphs}) is
extremely fast --- for 10 tetrahedra it takes a little over five minutes
on a 2.40GHz Pentium~4.  As such, there is no urgency to improve the
efficiency of this enumeration at the present time; a full list of
candidate face pairing graphs has been obtained for up to 13 tetrahedra,
well beyond the current census limits.
Table~\ref{tab-graphs} presents the number of connected 4-valent
multigraphs in these lists.

\begin{table}[htb]
    \caption{The number of connected 4-valent multigraphs for $\leq 13$
        tetrahedra}
    \label{tab-graphs}
    \begin{center}
        \begin{tabular}{|c|r|}
        \hline
        \bf \# Tet. & \bf Graphs \\
        \hline
        1 & 1 \\ 2 & 2 \\ 3 & 4 \\ 4 & 10 \\ 5 & 28 \\
        \hline
        \end{tabular}
        \qquad
        \begin{tabular}{|c|r|}
        \hline
        \bf \# Tet. & \bf Graphs \\
        \hline
        6 & 97 \\ 7 & 359 \\ 8 & 1\,635 \\ 9 & 8\,296 \\ 10 & 48\,432 \\
        \hline
        \end{tabular}
        \qquad
        \begin{tabular}{|c|r|}
        \hline
        \bf \# Tet. & {} \hfill \bf Graphs \hfill {} \\
        \hline
        11 & 316\,520 \\ 12 & 2\,305\,104 \\ 13 & 18\,428\,254 \\
        \hline
        \multicolumn{2}{c}{~} \\
        \multicolumn{2}{c}{~}
    \end{tabular} \end{center}
\end{table}

On the other hand, the processing of face pairing graphs
(step~\ref{e-processgraphs}) is extremely slow, consuming virtually the
entire running time of the generation algorithm.  For the 10-tetrahedron
non-orientable census this processing required approximately
$3 \frac23$ years of CPU time, though this was distributed amongst many
machines in parallel with a final cost of only
$1 \frac12$ months by the wall clock.

The division of tasks in Algorithm~\ref{a-generation}
lends itself to two avenues for improvement, both
of which are explored in this paper.
\begin{itemize}
    \item We can prove constraints that must be satisfied by the face
    pairing graph of a closed minimal {\ppirr} 3-manifold triangulation.
    Graphs that do not satisfy these constraints can be discarded
    after step~\ref{e-enumgraphs}, and do not need to be processed.

    Several results of this type are proven in \cite{burton-facegraphs},
    and we prove more in Section~\ref{s-graphs}.  Since the enumeration
    of graphs takes negligible time, if we can avoid processing $k\%$ of
    potential graphs in this way then we can expect to save
    roughly $k\%$ of the overall running time.

    \item We can improve the processing itself, finding ways to prune
    the search for possible triangulations in step~\ref{e-processgraphs}
    of the algorithm.  In Section~\ref{s-ufind} we accomplish this by
    using a modified union-find algorithm that tracks vertex and edge
    links.  As a consequence we are able to decrease running times by
    orders of magnitude, as seen in the experimental results of
    Section~\ref{s-ufind-timing}.
    %
\end{itemize}

%% file: prelim.tex
\section{Preliminary Results} \label{s-prelim}

Before moving on to specific algorithmic improvements, we present a
number of useful properties of minimal {\ppirr} triangulations.
These are simple structural properties that are easily tested by
computer, and many will be called upon in later sections of this paper.
We begin with a number of results from the literature, and then present some
new variants on existing results.

Our first property constrains the number of vertices in a minimal
{\ppirr} triangulation.  The following lemma was proven in the orientable
case by Jaco and Rubinstein \cite{0-efficiency}, and an equivalent
result involving special spines was proven by
Martelli and Petronio \cite{italian-decomp} for both orientable and
non-orientable manifolds.

\begin{lemma} \label{l-onevertex}
    Let $M$ be a closed {\ppirr} 3-manifold that is not
    $\sss$, $\rps$ or $L_{3,1}$.
    Then every minimal triangulation of $M$ contains precisely one vertex.
\end{lemma}

It will prove useful in later sections to count edges as well as vertices.
A quick Euler characteristic calculation yields the following immediate
consequence.

\begin{corollary} \label{c-edgecount}
    Let $T$ be a minimal triangulation of a
    closed {\ppirr} 3-manifold that is not
    $\sss$, $\rps$ or $L_{3,1}$, and
    suppose that $T$ is formed from $n$ tetrahedra.  Then
    $T$ contains precisely $n+1$ edges.
\end{corollary}

The following structural results are proven for both
orientable and non-orientable triangulations by this author in
\cite{burton-facegraphs}.
When restricted to the orientable case, a number of similar results
have been proven by other authors; see \cite{burton-facegraphs}
again for full details.

\begin{lemma} \label{l-edgedeg}
    Let $T$ be a closed minimal {\ppirr} triangulation containing
    $\geq 3$ tetrahedra.
    Then no edge of $T$ has degree one or two.
    Moreover, if an edge of $T$ has degree three then that edge
    does not lie on three distinct tetrahedra (instead, it
    meets some tetrahedron more than once).
\end{lemma}

\begin{figure}[htb]
    \centerline{\includegraphics[scale=0.7]{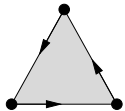}}
    \caption{All three edges of a face identified with each other}
    \label{fig-l31}
\end{figure}

\begin{lemma} \label{l-l31}
    Let $T$ be a closed minimal {\ppirr} triangulation containing
    $\geq 3$ tetrahedra.
    Then no face of $T$ has all three of its edges identified with each
    other in the same direction around the face,
    as illustrated in Figure~\ref{fig-l31}.
\end{lemma}

\begin{figure}[htb]
    \centerline{\includegraphics[scale=0.7]{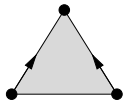}}
    \caption{One face with its edges identified to form a cone}
    \label{fig-cone}
\end{figure}

\begin{lemma} \label{l-cone}
    Let $T$ be a closed minimal {\ppirr} triangulation containing
    $\geq 3$ tetrahedra.
    Then no face of $T$ has two of its edges identified to form a cone,
    as illustrated in Figure~\ref{fig-cone}.
\end{lemma}

The next result regarding two-face cones is proven in a more general
form by the author in \cite{burton-facegraphs}.  Here we only require a
simpler variant in which the interior edges of the cone are distinct.

\begin{figure}[htb]
    \centerline{\includegraphics[scale=0.7]{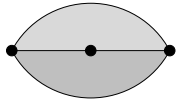}}
    \caption{Two faces joined together to form a cone}
    \label{fig-cone2}
\end{figure}

\begin{lemma} \label{l-cone2}
    Let $T$ be a triangulation of a closed 3-manifold containing
    $\geq 3$ tetrahedra.  Suppose that two distinct faces of $T$
    are joined together along their edges to form a cone, as illustrated
    in Figure~\ref{fig-cone2}.  Moreover, suppose that the two internal
    edges of this cone are distinct (i.e., they are not identified in
    the overall triangulation).
    Then either $T$ is non-minimal, or the 3-manifold that $T$
    represents is not {\ppirr}.
\end{lemma}

Having cited results regarding one-face cones and two-face cones, we
proceed to prove a similar result regarding three-face cones.  A
restriction of this result to orientable triangulations is proven by
Martelli and Petronio in \cite{italian9}.

\begin{figure}[htb]
    \centerline{\includegraphics[scale=0.5]{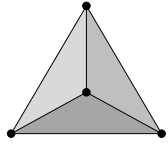}}
    \caption{Three faces joined together to form a cone}
    \label{fig-cone3}
\end{figure}

\begin{figure}[htb]
    \psfrag{Tip}{{\small Tip of cone}}
    \centerline{\includegraphics[scale=0.6]{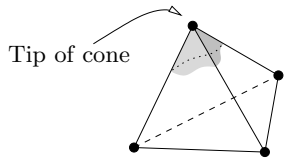}}
    \caption{A three-face cone surrounding a single tetrahedron vertex}
    \label{fig-cone3tet}
\end{figure}

\begin{lemma} \label{l-cone3}
    Let $T$ be a closed minimal {\ppirr} triangulation containing
    $\geq 3$ tetrahedra.
    Suppose that three faces of $T$ are joined together along their
    edges to form a cone, as illustrated in Figure~\ref{fig-cone3}.
    Moreover, suppose that the three internal edges of this cone are
    distinct (i.e., they are not identified in the overall triangulation).

    Then the three faces of this cone are distinct faces of the
    triangulation, and all three faces belong to a single tetrahedron.
    More specifically, this cone is the cone surrounding some vertex of
    the tetrahedron, as illustrated in Figure~\ref{fig-cone3tet}.
\end{lemma}

\begin{proof}
    It is simple to see that the three faces of the cone are distinct.
    A quick run through the six possible ways in which two faces could
    be identified shows that, in each case, either two internal edges
    are also identified (contradicting the conditions of the lemma)
    or some internal edge is identified with itself in reverse
    (giving a structure that is not a 3-manifold triangulation).

    Suppose then that the three-face cone does not surround
    some vertex of a single tetrahedron.  Let $n$ be the number of
    tetrahedra in triangulation $T$.
    Our strategy is to expand the triangulation in the
    region of the cone, and then to compress it along an edge to obtain a
    new triangulation of the same 3-manifold that uses fewer than $n$
    tetrahedra.  This mirrors the procedure used by Martelli and Petronio
    in the orientable case.  We take each step carefully in order to
    avoid corner cases that might behave unexpectedly.

    Since the three interior edges of the cone are distinct,
    the interior of the cone is embedded within the triangulation.
    This allows us to thicken the interior of the cone
    as follows.  We split the centre of the cone into two separate
    vertices and pull these vertices apart.  Likewise, the original cone
    is pulled apart into two cones (an upper and a lower cone), joined along
    their three boundary edges.
    In the empty space that is created, we insert two tetrahedra.

    \begin{figure}[htb]
        \psfrag{U}{{\small $U$}} \psfrag{V}{{\small $V$}}
        \psfrag{W}{{\small $W$}}
        \psfrag{X1}{{\small $X_1$}} \psfrag{X2}{{\small $X_2$}}
        \psfrag{X3}{{\small $X_3$}}
        \psfrag{Upper}{{\small Upper tetrahedron}}
        \psfrag{Lower}{{\small Lower tetrahedron}}
        \centerline{\includegraphics[scale=0.6]{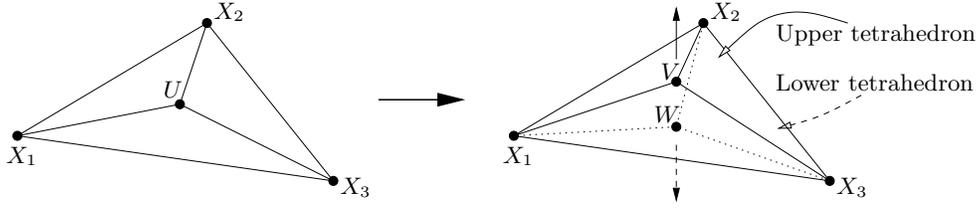}}
        \caption{Expanding a three-face cone into a region bounding two
            tetrahedra}
        \label{fig-cone3thicken}
    \end{figure}

    This expansion is illustrated in Figure~\ref{fig-cone3thicken},
    where the original vertex $U$ is split into two vertices $V$ and $W$.
    The two new tetrahedra are separated by the internal face
    $\Delta X_1X_2X_3$,
    with the upper tetrahedron containing vertices $V$, $X_1$, $X_2$ and
    $X_3$, and the lower tetrahedron containing vertices
    $W$, $X_1$, $X_2$ and $X_3$.

    The result is a new triangulation $T'$ that represents the same
    3-manifold as $T$ using precisely $n+2$ tetrahedra.
    More importantly, the two new vertices $V$ and $W$ are distinct.  Our
    aim now is to find an edge joining two distinct vertices
    that is suitable for compression.

    Since the original cone does not surround the vertex of a single
    tetrahedron either above or below, it follows that the new upper faces
    $\Delta X_1X_2V$, $\Delta X_2X_3V$ and $\Delta X_3X_1V$ collectively
    meet at least three distinct tetrahedra, including the new tetrahedron
    $X_1X_2X_3V$ (the only remaining possibilities for $\leq 2$
    distinct tetrahedra give a vertex link for $V$ that is not a 2-sphere).
    Likewise, the lower faces
    $\Delta X_1X_2W$, $\Delta X_2X_3W$ and $\Delta X_3X_1W$ collectively
    meet at least three distinct tetrahedra.
    Therefore at least two of the edges
    $\edge{X_1V}$, $\edge{X_2V}$ and $\edge{X_3V}$
    meet $\geq 3$ distinct tetrahedra each, and at least two of the
    edges $\edge{X_1W}$, $\edge{X_2W}$ and $\edge{X_3W}$
    meet $\geq 3$ distinct tetrahedra each.

    It follows that for some $i \in \{1,2,3\}$, each of the edges
    $\edge{X_iV}$ and $\edge{X_iW}$ meets $\geq 3$ distinct
    tetrahedra.  Since vertices $V$ and $W$ are distinct, one of these
    two edges must join two distinct vertices in the triangulation.
    Call this edge $e$.

    Our new aim is to compress the expanded triangulation $T'$ along
    edge $e$, producing a new triangulation $T''$ with $<n$ tetrahedra.
    The compression is performed as follows.

    \begin{enumerate}
        \item We shrink edge $e$ to a single point.  Since $e$ joins two
        distinct vertices, this does not alter the underlying 3-manifold.
        However, every tetrahedron that contains $e$ changes its shape
        according to Figure~\ref{fig-edgecomp}.

        \begin{figure}[htb]
            \psfrag{e}{{\small $e$}}
            \centerline{\includegraphics[scale=0.6]{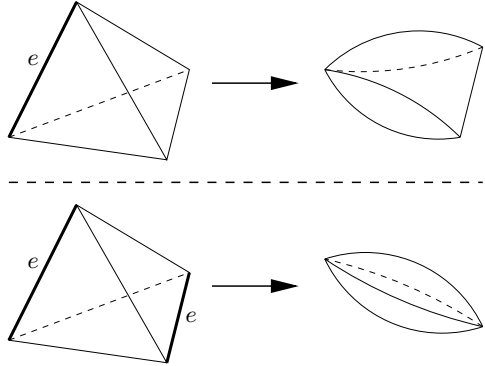}}
            \caption{Consequences of shrinking edge $e$ to a single point}
            \label{fig-edgecomp}
        \end{figure}

        The possible ways in which a tetrahedron may contain edge $e$
        one or more times are as follows.
        \begin{itemize}
            \item The tetrahedron contains edge $e$ only once.
            In this case the result is a triangular pillow with two
            bigon sides, as seen in the upper diagram of
            Figure~\ref{fig-edgecomp}.

            \item The tetrahedron contains edge $e$ twice, where $e$
            appears as two opposite edges of the tetrahedron.
            The result here is a football with four bigon faces, as seen
            in the lower diagram of Figure~\ref{fig-edgecomp}.

            \item The tetrahedron contains edge $e$ more than once, and
            $e$ appears twice on a single face.  Since the endpoints of
            $e$ are distinct vertices, the face must have its edges $e$
            identified to form a cone as illustrated earlier in
            Figure~\ref{fig-cone}.  This means that, unless the face is
            $\Delta X_1X_2X_3$, there is a corresponding face in $T$
            that also forms a cone, in contradiction to Lemma~\ref{l-cone}.

            Thus edge $e$ belongs to face $\Delta X_1X_2X_3$ more than
            once, forming a one-face cone in $T'$ but not in $T$.
            Without loss of generality, let directed edge
            $\vec{e} = \dedge{X_1X_2} = \dedge{X_3X_2}$ in the expanded
            triangulation $T'$.  This implies that
            $\dedge{X_1X_2} = \dedge{X_3X_2}$ in the original
            triangulation $T$.

            To avoid a two-face cone that
            contradicts Lemma~\ref{l-cone2}, these edges must also be
            identified with $\edge{UX_2}$ in some direction.
            To avoid a one-face cone that contradicts Lemma~\ref{l-cone},
            this direction must be
            $\dedge{X_1X_2}=\dedge{X_3X_2}=\dedge{X_2U}$.
            Moving back to the expanded triangulation $T'$, the old edge
            $\dedge{X_2U}$ is pulled apart into two distinct edges
            $\dedge{X_2V}$ and $\dedge{X_2W}$.  However, one of these
            identifications with $\vec{e}$ will be preserved; without loss of
            generality, suppose that $\vec{e}$ is identified with the new
            directed edge $\dedge{X_2V}$.  This gives
            $\vec{e}=\dedge{X_1X_2}=\dedge{X_2V}$ in $T'$, contradicting the
            fact that the endpoints of $e$ are distinct.
        \end{itemize}

        \item We eliminate each football with four bigon faces, using a
        combination of the following two operations. \label{en-flatfootball}
        The order in which
        these operations are performed is important, and will be
        discussed shortly.  In the meantime however, we describe each
        operation in detail.

        \begin{enumerate}[(i)]
            \item {\em Flattening bigon faces to edges.}
            This procedure is illustrated in Figure~\ref{fig-flatbigon}.

            \begin{figure}[htb]
                \centerline{\includegraphics[scale=0.6]{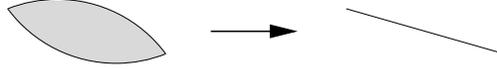}}
                \caption{Flattening a bigon face to a single edge}
                \label{fig-flatbigon}
            \end{figure}

            Flattening a bigon face has no effect upon the
            underlying 3-manifold unless the two edges of the bigon are
            identified.  In this case, the bigon forms either a 2-sphere
            or a projective plane.

            \begin{itemize}
                \item If the bigon forms a 2-sphere, flattening the
                bigon has the effect of cutting the manifold along this
                sphere and filling the resulting boundaries with solid
                balls.  Since the original manifold is {\ppirr},
                the resulting manifold is a disconnected union of the
                original manifold and a 3-sphere.  In this case
                we simply toss away the extra 3-sphere (possibly reducing
                the overall number of tetrahedra) and keep the structure
                that corresponds to the original manifold.

                \item If the bigon forms a projective plane then this plane
                has just one vertex, one edge and one face, and is
                therefore embedded in the underlying 3-manifold.
                The plane cannot be two-sided, since the underlying
                manifold is {\ppirr}.  Therefore the manifold contains
                an embedded one-sided projective plane, and is thus
                of the form $M \# \rps$ for some manifold $M$.  From
                {\ppirrty} again it follows that $M$ is trivial and the
                original manifold is simply $\rps$.  This however can
                be triangulated with two tetrahedra \cite{burton-facegraphs},
                contradicting the requirement that the original
                triangulation $T$ be minimal.
            \end{itemize}

            \item {\em Flattening bigon pillows to bigon faces.}
            This procedure is illustrated in Figure~\ref{fig-flatbigonpillow}.

            \begin{figure}[htb]
                \centerline{\includegraphics[scale=0.6]{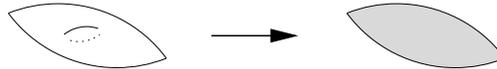}}
                \caption{Flattening a bigon pillow to a single bigon face}
                \label{fig-flatbigonpillow}
            \end{figure}

            Flattening a bigon pillow cannot change the underlying
            3-manifold unless the two faces of the pillow are identified.
            If this is the case then the bigon pillow represents the entire
            3-manifold, since it forms a complete closed structure.  Running
            through the four possible ways in which the upper bigon face can
            be identified with the lower, we see that the 3-manifold must be
            either $\sss$ or $\rps$.  Both manifolds can be triangulated
            using $\leq 2$ tetrahedra \cite{burton-facegraphs}, in
            contradiction to the
            requirement that triangulation $T$ be minimal.
        \end{enumerate}

        Having established that neither operation changes the underlying
        3-manifold, we can eliminate all footballs using the following
        algorithm.

        \begin{itemize}
            \item The overall aim of the algorithm is to
            eliminate each football with four bigon faces by
            (i)~flattening two of the four bigon faces to edges, thereby
            reducing the football to a bigon pillow, and then
            (ii)~flattening this bigon pillow to a bigon face, removing
            the football completely.

            \item The algorithm consists of a single operation that is
            repeated over and over.  Each time we begin the operation,
            all remaining footballs have either three or four bigon
            faces --- this is certainly true at the very beginning of
            the algorithm, and it will be seen that this remains true
            each time the operation is applied.

            \item The operation is as follows.  Locate a bigon face $B$
            that belongs to one of the remaining footballs, where $B$ is
            not identified with some other face of the same football.

            We flatten this bigon $B$ to an edge.  This will reduce the
            number of faces of the chosen football by one.  If $B$ was
            identified with a bigon face of some different football,
            this other football will also have its face count reduced by one.
            Otherwise $B$ was identified with a face of a triangular
            pillow, which will thereby lose one of its bigon sides.

            At this point we potentially have one or two footballs with
            only two bigon faces, i.e., we potentially have one or two
            bigon pillows.  If this occurs, the bigon pillows are
            immediately removed by flattening them to faces.  All
            remaining footballs will therefore have either three or four
            bigon faces as required, and the operation can begin again.

            \item We continue until the operation can no longer be
            performed.  If any footballs still remain, they must be
            four-face footballs whose four bigon faces are identified in
            two pairs.
            However, such a football is impossible --- we could flatten
            the first bigon to an edge, reducing the football to a bigon
            pillow whose two remaining
            faces are identified, and this was seen to
            be impossible in step~(ii) above.
        \end{itemize}

        Once this algorithm has finished,
        the footballs will have disappeared completely.
        In addition, some of the triangular pillows with bigon sides might
        have had their sides flattened to edges.
        The possible new shapes that these triangular pillows might take
        can be seen in Figure~\ref{fig-alltripillows}.

        \begin{figure}[htb]
            \centerline{\includegraphics[scale=0.6]{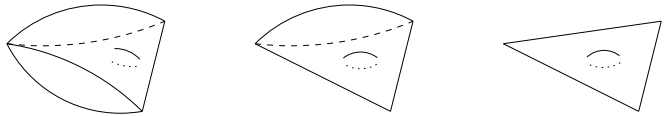}}
            \caption{Remaining triangular pillows with or
                without bigon sides}
            \label{fig-alltripillows}
        \end{figure}

        \item We flatten all remaining bigon faces to edges.  This
        simplifies the triangular pillows so that each pillow has just
        three edges and two triangular faces, as seen in the rightmost
        diagram of Figure~\ref{fig-alltripillows}.

        \item We flatten all triangular pillows to faces.  This
        procedure is illustrated in Figure~\ref{fig-flattripillow}.

        \begin{figure}[htb]
            \centerline{\includegraphics[scale=0.6]{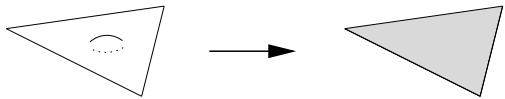}}
            \caption{Flattening a triangular pillow to a single face}
            \label{fig-flattripillow}
        \end{figure}

        The only way in which flattening a triangular pillow
        can change the underlying 3-manifold is if the upper and lower
        faces of the pillow are identified.  In such a case, the pillow
        must represent the entire 3-manifold (since it forms a closed
        structure), and a quick run through the six possible
        rotations and reflections of the triangle shows that this
        3-manifold is either $\sss$ or $L_{3,1}$.  Both manifolds can be
        triangulated using $\leq 2$ tetrahedra as seen in
        \cite{burton-facegraphs}, contradicting the claim that the
        original triangulation $T$ is minimal.
    \end{enumerate}

    At this point all of our unusual shapes have been flattened away,
    and we are left with only tetrahedra.  This gives us a new
    triangulation $T''$ of the original 3-manifold, where each tetrahedron
    of $T'$ that contained edge $e$ has been removed.  Since
    triangulation $T'$ contained precisely $n+2$ tetrahedra and edge $e$
    was chosen to belong to $\geq 3$ distinct tetrahedra, it follows
    that our final triangulation $T''$ contains at most $n-1$ tetrahedra.
    Therefore the original triangulation $T$ could not have been
    minimal, and the proof is complete.
\end{proof}

To close this section, we examine some consequences of the previous
results when applied to tori within a triangulation.

\begin{figure}[htb]
    \psfrag{a}{{\small $a$}} \psfrag{b}{{\small $b$}} \psfrag{c}{{\small $c$}}
    \psfrag{A}{{\small $A$}} \psfrag{B}{{\small $B$}}
    \psfrag{C}{{\small $C$}} \psfrag{D}{{\small $D$}}
    \centerline{\includegraphics[scale=0.6]{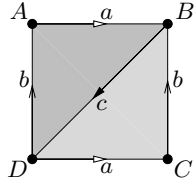}}
    \caption{Two faces joined together to form a torus}
    \label{fig-torus}
\end{figure}

\begin{corollary} \label{c-torusdistinct}
    Let $T$ be a closed minimal {\ppirr} triangulation containing
    $\geq 3$ tetrahedra.  Suppose that two distinct faces of $T$ are
    joined along their edges to form a torus, as illustrated in
    Figure~\ref{fig-torus}.  Then all three edges of this torus are
    distinct, i.e., no two of these edges are identified in the overall
    triangulation.
\end{corollary}

\begin{proof}
    It is clear that all three edges of the torus cannot be identified
    together, since this would break either Lemma~\ref{l-l31} or~\ref{l-cone}.
    Suppose then, without loss of generality, that edges $a$ and $b$ are
    identified.

    If directed edges $a$ and $-b$ are identified then each face
    has two edges identified to form a cone, in
    contradiction to Lemma~\ref{l-cone}.  It follows that directed
    edges $a$ and $b$ are identified.  This however creates a two-face
    cone centred about vertex $D$, with internal edges $a=b$ and $c$
    (where directed edge $a=b$ points away from the centre of the cone,
    and directed edge $c$ points towards the centre of the cone).

    Since all three edges of the torus are distinct, it follows that
    the internal edges of this cone are likewise distinct.  This contradicts
    Lemma~\ref{l-cone2} and concludes the proof.
\end{proof}

\begin{corollary} \label{c-toruslayer}
    Let $T$ be a closed minimal {\ppirr} triangulation containing
    $\geq 3$ tetrahedra.  Suppose that two distinct faces of $T$ are
    joined along their edges to form a torus, as illustrated previously
    in Figure~\ref{fig-torus}.  Moreover, suppose that some other face
    of $T$ contains at least two of the edges $a,b,c$.  Then all three
    faces belong to a single tetrahedron.
\end{corollary}

\begin{proof}
    Denote the third face by $F$.
    Without loss of generality, suppose that $F$ contains directed
    edges $a$ and $b$.  Noting from Corollary~\ref{c-torusdistinct} that
    $a$ and $b$ are distinct, Figure~\ref{fig-toruslayercases} shows
    the possible ways in which these edges may be oriented within face $F$.

    \begin{figure}[htb]
        \psfrag{a}{{\small $a$}} \psfrag{b}{{\small $b$}}
        \psfrag{(i)}{{\small (i)}} \psfrag{(ii)}{{\small (ii)}}
        \psfrag{(iii)}{{\small (iii)}} \psfrag{(iv)}{{\small (iv)}}
        \centerline{\includegraphics[scale=0.6]{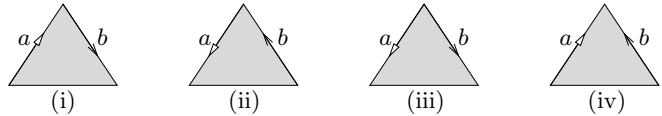}}
        \caption{Possible cases for directed edges $a$ and $b$ within
            the third face $F$}
        \label{fig-toruslayercases}
    \end{figure}

    In case~(i), face $F$ forms a two-face cone with face \polytri{DCB}
    from the torus (with vertex $C$ at the centre of the cone).
    Similarly, in case~(ii), face $F$ forms a two-face cone with face
    \polytri{BAD}.  Lemma~\ref{l-cone2} provides a contradiction in
    both cases.

    For case~(iii), we obtain a cone from all three faces.  Vertex $D$
    sits at the centre of this cone, with directed edges $a$ and $b$
    running away from the centre and $c$ running towards the centre.
    Since $a$, $b$ and $c$ are all distinct from
    Corollary~\ref{c-torusdistinct}, we see from Lemma~\ref{l-cone3}
    that all three faces belong to a single tetrahedron.
    Case~(iv) is handled likewise.
\end{proof}

%% file: graphs.tex
\section{Face Pairing Graphs} \label{s-graphs}

Recall from Definition~\ref{d-facegraph} that the {\em face pairing graph}
of a triangulation is a 4-valent graph describing which tetrahedron
faces are identified with which.  In this section we use an analysis of
face pairing graphs to obtain our first round of improvements to the
census algorithm.

The earlier paper \cite{burton-facegraphs} likewise uses face pairing
graphs to improve the census algorithm; many of the results proven here
are generalisations of these earlier theorems.  Some core results
of \cite{burton-facegraphs} are reviewed in Section~\ref{s-graphs-prev},
where we also introduce the recurring concepts of chains and layered solid
tori.

Following this, Section~\ref{s-graphs-elim}
builds on this work to obtain new results and algorithmic improvements.
In this section we obtain a series of properties
that must be satisfied by the face pairing graph of a closed minimal
{\ppirr} triangulation.  Using these results, we can improve
Algorithm~\ref{a-generation} by identifying graphs that do not
satisfy these properties and eliminating them
after step~\ref{e-enumgraphs} of the
algorithm (enumeration of graphs), without carrying them through
to step~\ref{e-processgraphs} (processing of graphs).

Section~\ref{s-graphs-data} gives an indication of the usefulness of
these results.  In particular, it examines all possible face pairing
graphs on $\leq 10$ vertices and identifies how many of these graphs can
be eliminated using the results of Section~\ref{s-graphs-elim}, thereby
offering a rough estimate of the processing time that can be saved.


\subsection{Previous Results} \label{s-graphs-prev}

Before summarising the key results of \cite{burton-facegraphs}, some
additional terminology is required.  Most of the major results of
\cite{burton-facegraphs} (and all of the graph-related results of this
paper) involve subgraphs called {\em chains}.

\begin{defn}[Chain]
    A {\em chain} of length $k$ is a sequence of $k$ double edges joined
    end-to-end in a linear fashion, as illustrated in Figure~\ref{fig-chain}.
    A chain of length $k$ must link together $k+1$ distinct vertices
    (i.e., no vertex may be repeated).

    \begin{figure}[htb]
        \centerline{\includegraphics[scale=0.7]{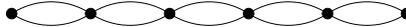}}
        \caption{A chain of length five}
        \label{fig-chain}
    \end{figure}

    A {\em one-ended chain} of length $k$ is a chain of length $k$ with
    a loop added to one end, as illustrated in Figure~\ref{fig-chainoneend}.
    Likewise, a {\em double-ended chain} of length $k$ is a chain of
    length $k$ with loops added to both ends.

    \begin{figure}[htb]
        \centerline{\includegraphics[scale=0.7]{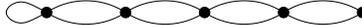}}
        \caption{A one-ended chain of length four}
        \label{fig-chainoneend}
    \end{figure}

    Wherever this paper refers to a one-ended or double-ended chain,
    chains of length zero are also included.  In particular, any
    statement involving one-ended chains is assumed to include
    the zero-length case of a single loop.
\end{defn}

A key property of chains is that every vertex has degree four except for
the two vertices at each end.  This means that, within a 4-valent face pairing
graph, a chain is almost entirely self-contained --- it may only join to
the remainder of the graph at these two end vertices.

The value of this is seen when we translate to the language of
triangulations.  Here a chain of length $k$ becomes a
sequence of $k+1$ tetrahedra joined together, with almost all of these
tetrahedra meeting other tetrahedra from the chain along all four faces.
The only tetrahedra with any spare faces for connecting this structure to
the rest of the triangulation are the two tetrahedra at either end of
the sequence.
This constrains the possibilities enough for us to draw interesting conclusions
about the rotations and reflections used when identifying tetrahedron
faces together.

One-ended chains are even more constrained, since every vertex but
one has degree four.
This allows the corresponding tetrahedron-based
structures to be identified even more precisely, as seen in
\cite{burton-facegraphs} and reproduced in Theorem~\ref{t-chain-lst} below.
To understand these structures, we must introduce the concept of
layering.

\begin{defn}[Layering]
    Consider a triangulation with boundary $B$, and let $e$ be an edge
    of this boundary that lies between two distinct boundary faces.
    To {\em layer on edge $e$} is to perform the following operation.

    We obtain a new tetrahedron $\Delta$, and lay two adjacent faces
    of $\Delta$ directly onto the two boundary faces surrounding $e$.
    These two original faces surrounding $e$ become internal faces of the
    triangulation, and the remaining two faces of $\Delta$ become new
    boundary faces.  Likewise, the original boundary edge $e$ becomes
    internal, making way for a new boundary edge $f$ on the far side of
    tetrahedron $\Delta$.  This procedure is illustrated in
    Figure~\ref{fig-layering}.

    \begin{figure}[htb]
        \psfrag{e}{{\small $e$}}
        \psfrag{f}{{\small $f$}}
        \psfrag{T}{{\small $\Delta$}}
        \psfrag{B}{{\small $B$}}
        \centerline{\includegraphics[scale=0.7]{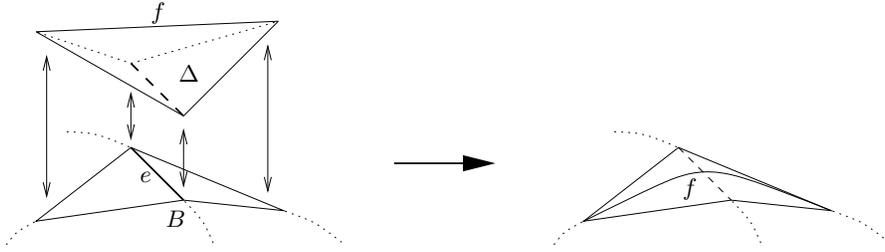}}
        \caption{Performing a layering}
        \label{fig-layering}
    \end{figure}
\end{defn}

It should be observed that performing a layering does not alter the
underlying 3-manifold of a triangulation.  What it does change is the
curves made by the boundary edges of the triangulation, since in general
the new boundary edge $f$ does not represent the same curve as the original
boundary edge $e$.

One of the more frequent examples of layering that can be found in a
census of triangulations is in the construction of {\em layered solid tori}.
These well-structured triangulations of solid tori have been discussed
in depth by Jaco and Rubinstein \cite{0-efficiency,layeredlensspaces},
and similar structures involving special spines have been described by
Matveev \cite{matveev6}.  In the context of a census of 3-manifold
triangulations, they are parameterised in detail in \cite{burton-nor7}.

We do not require such detail here, and so a full definition of a
layered solid torus is omitted --- the reader is referred to the
aforementioned references for details.  For this paper the following
summary is sufficient.

\begin{defn}[Layered Solid Torus]
    A {\em standard layered solid torus} is a triangulation of the solid
    torus obtained through successive layerings upon the boundary of $T_1$,
    where $T_1$ is the one-tetrahedron triangulation of the solid torus.

    There is no limit upon the number of layerings --- there may be a
    large number of layerings or there may be none at all.  There is
    also no constraint upon which particular boundary edges are used for
    performing the layerings.

    As a result, an infinite variety of different standard layered solid
    tori can be constructed.  Each of these solid tori has two boundary
    faces, arranged to form a torus as illustrated in
    Figure~\ref{fig-lst}.

    \begin{figure}[htb]
    \psfrag{Torus}{{\small \twolines{c}{Solid}{torus}}}
    \centerline{\includegraphics[scale=0.7]{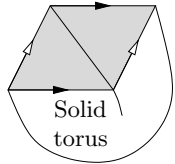}}
    \caption{The boundary of a layered solid torus}
    \label{fig-lst}
    \end{figure}
\end{defn}

Equipped with this terminology, we are now in a position to recount the
core results proven by the author in \cite{burton-facegraphs}.
The first of these results describes the tetrahedron-based structures
corresponding to one-ended chains.  This theorem forms a starting point
for many later results, including all of the results of
Section~\ref{s-graphs-elim} in this paper.

\begin{theorem} \label{t-chain-lst}
    Let $G$ be the face pairing graph of a closed minimal {\ppirr}
    triangulation containing $\geq 3$ tetrahedra.  Suppose that $G$
    contains a one-ended chain.  Then the tetrahedra corresponding to
    the vertices of this one-ended chain form a standard layered solid
    torus within the triangulation.
\end{theorem}

Recall that every vertex of a one-ended chain has degree four except for
a single end vertex, which has degree two.  This allows a one-ended
chain to be attached to the remainder of a face pairing graph by two
edges; in the language of triangulations, this corresponds to attaching
a layered solid torus to the remainder of the triangulation along its
two boundary faces.

In Section~3.2 of \cite{burton-facegraphs}, the author
proves a series of constraints on the face pairing graphs of minimal
triangulations.  These constraints are bundled together into the single
theorem below.

\begin{theorem} \label{t-badgraphsold}
    Let $G$ be the face pairing graph of a closed minimal {\ppirr}
    triangulation containing $\geq 3$ tetrahedra.  Then $G$ does not
    contain any of the following subgraphs:
    \begin{enumerate}[(i)]
        \item A triple edge; \label{en-triple}
        \item A triangle with a double edge and a one-ended chain
        connected to the opposite vertex; \label{en-dblhandle}
        \item A single edge with a one-ended chain connected to each
        vertex, where this single edge does not form part of a double
        edge. \label{en-brokenchain}
    \end{enumerate}
\end{theorem}

\begin{figure}[htb]
    \psfrag{Subgraph (i)}{{\small Subgraph (\ref{en-triple})}}
    \psfrag{Subgraph (ii)}{{\small Subgraph (\ref{en-dblhandle})}}
    \psfrag{Subgraph (iii)}{{\small Subgraph (\ref{en-brokenchain})}}
    \centerline{\includegraphics[scale=0.7]{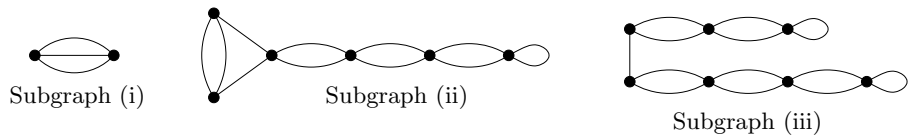}}
    \caption{Examples of prohibited subgraphs from earlier results}
    \label{fig-badgraphsold}
\end{figure}

Examples of these three prohibited types of subgraph are illustrated in
Figure~\ref{fig-badgraphsold}.
It will be seen in Section~\ref{s-graphs-elim} that several of the new
results proven in this paper can be viewed as generalisations of the
various components of the theorem above.

\subsection{Elimination of Graphs} \label{s-graphs-elim}

As indicated at the end of Section~\ref{s-algm}, one method of improving
the census algorithm involves proving constraints that must be satisfied
by the face pairing graphs of closed minimal {\ppirr} triangulations.
Results of this type allow us to identify graphs that can be discarded
immediately without any processing.

Theorem~\ref{t-badgraphsold} illustrates earlier results of this type
from \cite{burton-facegraphs}.  In this section we generalise these
earlier results, proving additional constraints in
Theorems~\ref{t-straybigon}, \ref{t-square} and~\ref{t-mountains}.
In the following section we present empirical data that measures how useful
these results might be in a real census.

It should be noted that in the orientable case, Martelli was
previously aware of simpler variants of some of these
constraints.\footnote{This was noted in a private communication with
the author in November, 2003.}
The results presented here are more general, and apply to
both orientable and non-orientable triangulations.

\begin{theorem} \label{t-straybigon}
    Let $G$ be the face pairing graph of a closed minimal {\ppirr}
    triangulation containing $\geq 3$ tetrahedra.  Suppose that $G$
    contain a one-ended chain, and suppose that some edge of $G$ joins
    this chain to another double edge as illustrated in
    Figure~\ref{fig-straybigon}.
    Then one of the following statements must be true:
    \begin{enumerate}[(i)]
        \item $G$ contains a longer one-ended chain, incorporating the
        original one-ended chain, edge and double edge previously
        discussed; \label{en-longerchain}
        \item The original one-ended chain is joined to both ends of a
        new chain of length two, where the original double edge
        belongs to this new chain; \label{en-doubledouble}
        \item There is some vertex $V$ not belonging to
        either the original chain or the double edge, where $V$ is
        joined by edges to the end of the original chain and to
        both ends of the double edge. \label{en-threesticks}
    \end{enumerate}

    \begin{figure}[htb]
        \centerline{\includegraphics[scale=0.7]{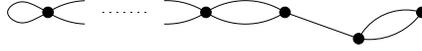}}
        \caption{Initial conditions for Theorem~\ref{t-straybigon}}
        \label{fig-straybigon}
    \end{figure}

    These three possibilities are illustrated in
    Figure~\ref{fig-straybigon-ok}, where the original one-ended chain,
    edge and double edge are drawn in solid lines, and any new vertices
    or edges are drawn in dashes.

    \begin{figure}[htb]
        \psfrag{Option (i)}{{\small Option (i)}}
        \psfrag{Option (ii)}{{\small Option (ii)}}
        \psfrag{Option (iii)}{{\small Option (iii)}}
        \psfrag{V}{{\small $V$}}
        \centerline{\includegraphics[scale=0.7]{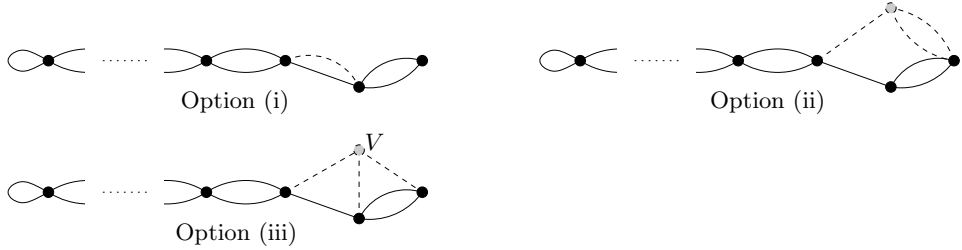}}
        \caption{Possible resolutions for Theorem~\ref{t-straybigon}}
        \label{fig-straybigon-ok}
    \end{figure}
\end{theorem}

\noindent
Before embarking on a proof, it is worth noting the following points
about Theorem~\ref{t-straybigon}.
\begin{itemize}
    \item Whilst a little unwieldy, all three options presented in
    the theorem statement are necessary.  Graphs corresponding to
    options~(\ref{en-longerchain}) and~(\ref{en-threesticks}) appear in
    very early levels of the orientable census; illustrated examples can
    be found in the six-tetrahedron census results of Matveev \cite{matveev6}.
    Option~(\ref{en-doubledouble}) is only required in the non-orientable
    case, where an 11-tetrahedron example has been found that
    necessitates it.
    \item This theorem generalises some of the earlier results of
    \cite{burton-facegraphs}.
    In particular, this new theorem covers both
    cases~(\ref{en-dblhandle}) and~(\ref{en-brokenchain}) of the earlier
    Theorem~\ref{t-badgraphsold}.
\end{itemize}

{ 
\renewcommand{\prooflabel}{Proof of Theorem~\ref{t-straybigon}}
\begin{proof}
    Suppose the face pairing graph $G$ contains a one-ended chain
    attached to some other double edge, as described in the theorem statement.
    Assign the label $V_1$ to the vertex at the non-loop end of the chain,
    and assign the labels $V_2$ and $V_3$ to the endpoints of the
    other double edge, where vertices $V_1$ and $V_2$ are joined.
    This is illustrated in the leftmost diagram of
    Figure~\ref{fig-straybigon-labels}.

    \begin{figure}[htb]
        \psfrag{V1}{{\small $V_1$}} \psfrag{V2}{{\small $V_2$}}
        \psfrag{V3}{{\small $V_3$}} \psfrag{V4}{{\small $V_4$}}
        \centerline{\includegraphics[scale=0.7]{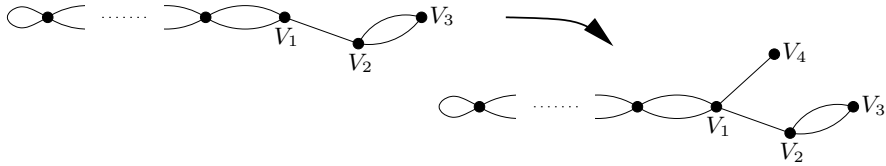}}
        \caption{Identifying graph vertices in Theorem~\ref{t-straybigon}}
        \label{fig-straybigon-labels}
    \end{figure}

    Since vertex $V_1$ only has degree three in this diagram, there is
    an edge meeting $V_1$ that is not yet accounted for.  This edge
    cannot join any other vertex of the chain, since these other
    vertices already have full degree four.  Likewise, it cannot be a
    loop joining $V_1$ to itself since this would raise the degree of
    $V_1$ to five.  If the edge runs from $V_1$ to $V_2$ then we have
    option~(\ref{en-longerchain}) of the theorem statement, and if it
    runs from $V_1$ to $V_3$ then we have a situation that
    Theorem~\ref{t-badgraphsold} proves impossible.

    The only case remaining is where $V_1$ is joined to some new vertex
    not previously considered.  Label this new vertex $V_4$; this is
    illustrated in the rightmost diagram of Figure~\ref{fig-straybigon-labels}.

    We translate now to the language of triangulations.  Suppose that
    vertices $V_1$, $V_2$, $V_3$ and $V_4$ of the graph correspond to
    tetrahedra $\Delta_1$, $\Delta_2$, $\Delta_3$ and $\Delta_4$
    of the triangulation respectively.
    From the graph we can draw the following conclusions.
    \begin{itemize}
        \item Two faces of tetrahedron $\Delta_1$ form the boundary
        torus of a standard layered solid torus (this follows from applying
        Theorem~\ref{t-chain-lst} to the one-ended chain).
        \item One of these faces is identified with some face of
        $\Delta_2$; the other is identified with some face of $\Delta_4$.
        \item Two other faces of $\Delta_2$ are each identified with
        faces of $\Delta_3$.
    \end{itemize}

    Part of this configuration is illustrated in
    Figure~\ref{fig-straybigon-above}.  The layered solid torus sits
    beneath the figure, with its two boundary faces \polytri{ABD} and
    \polytri{BDC} drawn in bold.  Note that edges \polyedge{AB} and
    \polyedge{DC} are identified, and likewise edges \polyedge{AD} and
    \polyedge{BC} are identified.  Tetrahedron $\Delta_4$ is then attached to
    face \polytri{ABD}, and uses vertices $A$, $B$, $D$ and $E$.
    Likewise, $\Delta_2$ is attached to face \polytri{BDC} and uses
    vertices $B$, $C$, $D$ and $F$.  Through the symmetry of the torus
    we may assume that faces \polytri{BDF} and \polytri{BCF} are joined
    to tetrahedron $\Delta_3$; in the rightmost diagram $\Delta_3$ is placed
    upon face \polytri{BDF}, using vertices $B$, $D$, $F$ and $G$.

    \begin{figure}[htb]
        \psfrag{A}{{\small $A$}} \psfrag{B}{{\small $B$}}
        \psfrag{C}{{\small $C$}} \psfrag{D}{{\small $D$}}
        \psfrag{E}{{\small $E$}} \psfrag{F}{{\small $F$}}
        \psfrag{G}{{\small $G$}}
        \psfrag{Label1}{{\small Faces of $\Delta_1$}}
        \psfrag{Label2}{{\small Attaching $\Delta_2$ and $\Delta_4$
            to $\Delta_1$}}
        \psfrag{Label3}{{\small Attaching $\Delta_3$ to $\Delta_2$}}
        \centerline{\includegraphics[scale=0.65]{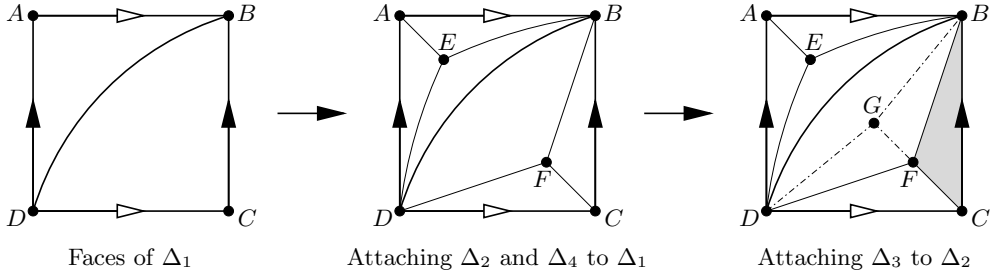}}
        \caption{Translating graph vertices to tetrahedra}
        \label{fig-straybigon-above}
    \end{figure}

    It remains then to identify \polytri{BCF} (shaded in the diagram)
    with one of the three
    remaining faces of $\Delta_3$.  Although there are 18 possibilities
    (three faces, with six rotations and reflections each), we can
    quickly cut this number down.

    Consider \polyedge{BC}; the final mapping must identify this edge with
    one of the edges of $\Delta_3$.  Noting that the distinct
    faces \polytri{ABD} and \polytri{BCD} form a torus,
    Corollary~\ref{c-torusdistinct} shows that \polyedge{BC} cannot be
    identified with \polyedge{BD}.  Corollary~\ref{c-toruslayer} then
    shows that \polyedge{BC} cannot be identified with \polyedge{BF} or
    \polyedge{FD}, since this would require \polytri{ABD}, \polytri{BDC}
    and \polytri{BDF} to belong to a single tetrahedron.  Likewise,
    Corollary~\ref{c-toruslayer} also shows that \polyedge{BC} cannot be
    identified with \polyedge{BG} or \polyedge{GD}, since otherwise
    \polytri{ABD}, \polytri{BDC} and \polytri{BDG} would belong to a
    single tetrahedron.

    The only remaining possibility is that \polyedge{BC} is identified
    with \polyedge{FG}.  Taking into account the two possible
    orientations of the edge and the two possible target faces
    \polytri{BFG} and \polytri{DFG}, this leaves us with the following
    four cases.  In each case, the order in which we write the vertices
    of the two identified faces indicates which particular rotation or
    reflection is used.\footnote{Specifically, a claim that
    faces \polytri{PQR} and \polytri{XYZ} are identified means that
    vertices $P$, $Q$ and $R$ are identified with $X$, $Y$ and
    $Z$ respectively.}
    \begin{itemize}
        \item \polytri{BCF} is identified with \polytri{FGB} (with
        directed edge \polydedge{BC} identified with \polydedge{FG}).

        In this case we have directed edges \polydedge{BF} and
        \polydedge{FB} identified, which is impossible.

        \item \polytri{BCF} is identified with \polytri{GFB} (with
        directed edge \polydedge{BC} identified with \polydedge{GF}).

        This forces the identification of directed edges
        $\dedge{GB} = \dedge{BF} = \dedge{FC}$, as illustrated in the
        leftmost diagram of Figure~\ref{fig-straybigon-cases}.
        Consider the three faces \polytri{ABD}, \polytri{DBG} and
        \polytri{FCD}.  Because of the aforementioned edge
        identifications, these three faces come together to form a
        three-face cone --- the three shaded arcs in the diagram mark
        the neighbourhood of the vertex at the centre of this cone.

        \begin{figure}[htb]
            \psfrag{A}{{\small $A$}} \psfrag{B}{{\small $B$}}
            \psfrag{C}{{\small $C$}} \psfrag{D}{{\small $D$}}
            \psfrag{E}{{\small $E$}} \psfrag{F}{{\small $F$}}
            \psfrag{G}{{\small $G$}}
            \psfrag{Label1}{{\small $\Delta \mathit{BCF}
                \leftrightarrow \Delta \mathit{GFB}$}}
            \psfrag{Label2}{{\small $\Delta \mathit{BCF}
                \leftrightarrow \Delta \mathit{FGD}$}}
            \psfrag{Label3}{{\small $\Delta \mathit{BCF}
                \leftrightarrow \Delta \mathit{GFD}$}}
            \centerline{\includegraphics[scale=0.65]{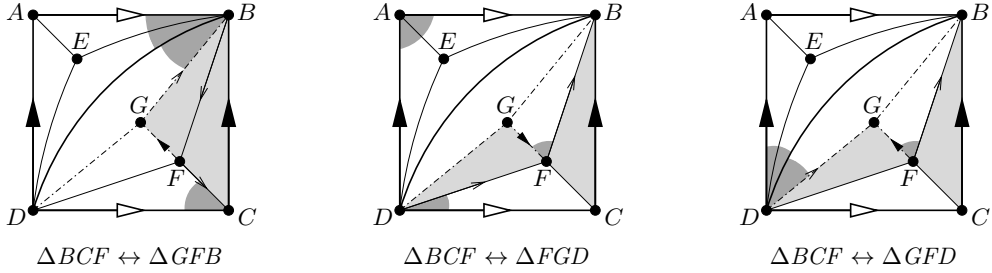}}
            \caption{Possibilities for the remaining join between
                $\Delta_2$ and $\Delta_3$}
            \label{fig-straybigon-cases}
        \end{figure}

        Our aim then is to use Lemma~\ref{l-cone3}.  First we observe
        that all three faces are distinct, since \polytri{ABD} is
        already an interior face, and identifying \polytri{DBG} with
        \polytri{FCD} would create a triple edge in the face pairing
        graph (which Theorem~\ref{t-badgraphsold} states is impossible).

        The three interior edges of the cone are \polyedge{AB},
        \polyedge{DB} and $\edge{GB}=\edge{BF}=\edge{FC}$.  We already
        know from Corollary~\ref{c-torusdistinct} that \polyedge{AB} and
        \polyedge{DB} are distinct; furthermore, since face
        \polytri{BCF} contains both \polyedge{GB} and the third edge of
        the original torus, Corollary~\ref{c-toruslayer} shows that
        \polyedge{GB} is distinct from both \polyedge{AB} and
        \polyedge{DB} (otherwise \polytri{BCF}, \polytri{ABD} and
        \polytri{BDC} would belong to a common tetrahedron, an
        impossible situation).

        Therefore we can apply Lemma~\ref{l-cone3} to the three-face
        cone described above.  It follows that \polytri{ABD},
        \polytri{DBG} and \polytri{FCD} belong to a single tetrahedron,
        which must be the upper-left tetrahedron $\Delta_4$.  Therefore
        vertices $V_2$ and $V_3$ each join to $V_4$ in the face pairing
        graph, and we have case~(\ref{en-threesticks}) from the theorem
        statement.

        \item \polytri{BCF} is identified with \polytri{FGD} (with
        directed edge \polydedge{BC} identified with \polydedge{FG}).

        As in the previous case, we seek to find and exploit a
        three-face cone.  This time we have directed edge
        identifications $\dedge{FB}=\dedge{DF}$ and
        $\dedge{GF}=\dedge{CB}=\dedge{DA}$, giving us a cone from faces
        \polytri{DAB}, \polytri{CDF} and \polytri{BFG}.  This is
        illustrated in the middle diagram of Figure~\ref{fig-straybigon-cases},
        with the shaded arcs once more representing a neighbourhood of the
        vertex at the centre of the cone.

        The three faces are distinct for the same reason as before
        (avoiding a triple edge in the face pairing graph), and again
        Corollaries~\ref{c-torusdistinct} and~\ref{c-toruslayer} show
        that the three internal edges of the cone are distinct.  Thus
        Lemma~\ref{l-cone3} shows that \polytri{DAB}, \polytri{CDF} and
        \polytri{BFG} belong to a common tetrahedron which must again
        be $\Delta_4$, giving us case~(\ref{en-threesticks}) from the
        theorem statement as before.

        \item \polytri{BCF} is identified with \polytri{GFD} (with
        directed edge \polydedge{BC} identified with \polydedge{GF}).

        In this final case, the relevant directed edge identifications
        are $\dedge{FB}=\dedge{DG}$ and $\dedge{FG}=\dedge{CB}=\dedge{DA}$.
        Again we find a three-face cone, this time using faces
        \polytri{ADB}, \polytri{BDG} and \polytri{BFG}.

        From here we diverge a little from previous cases.  This time
        the three faces are distinct because the alternative would be to
        identify two faces of $\Delta_3$, thereby adding a loop
        to $V_3$ in the face pairing graph and contradicting
        case~(\ref{en-brokenchain}) of Theorem~\ref{t-badgraphsold}.
        The three internal edges of the cone are distinct for the usual
        reasons, and so Lemma~\ref{l-cone3} shows that faces
        \polytri{ADB}, \polytri{BDG} and \polytri{BFG} all belong to a
        common tetrahedron.

        As before, the only possible common tetrahedron is the
        upper-left tetrahedron $\Delta_4$.  Since both
        \polytri{BDG} and \polytri{BFG} belong to $\Delta_3$, this adds
        a double edge from $V_3$ to $V_4$ in the face pairing graph,
        giving us case~(\ref{en-doubledouble}) from the theorem statement.
    \end{itemize}
    This concludes the list of possibilities, and the proof of
    Theorem~\ref{t-straybigon} is complete.
\end{proof}
} 

It should be noted that the face identifications in the final case
above give a non-orientable structure, justifying our earlier
claim that case~(\ref{en-doubledouble}) is only necessary in the
non-orientable case.  A little further investigation shows that the
manifold created by joining $\Delta_4$, $\Delta_3$, $\Delta_2$ and
the layered solid torus in this way is a bounded {\sfslong} of the form
$\sfs{\smash{\dischalfref}}{(p,q)}$,
where the orbifold $\dischalfref$ is the disc with
half-reflector, half-regular boundary.
See Section~\ref{s-census-nor} for further discussion of
manifolds of this type.

Our next two results extend cases~(\ref{en-triple}) and~(\ref{en-dblhandle})
of the old Theorem~\ref{t-badgraphsold}.  Recall that
case~(\ref{en-triple}) of Theorem~\ref{t-badgraphsold} outlaws a triple edge,
and case~(\ref{en-dblhandle}) can be viewed as outlawing a variant where one
of these three edges is replaced with a one-ended chain.  In
Theorems~\ref{t-square} and~\ref{t-mountains}, we likewise outlaw variants
in which two and then all three of the edges in a triple edge
are replaced with one-ended chains.

Since Theorems~\ref{t-square} and~\ref{t-mountains} are
similar in structure, we begin by stating the two theorems
without proof.  Following this we present and prove
Lemma~\ref{l-emptysquare}, which contains logic common to both theorems.
Finally we return to Theorems~\ref{t-square} and~\ref{t-mountains},
using this new lemma to prove each in turn.

\begin{theorem} \label{t-square}
    Let $G$ be the face pairing graph of a closed minimal {\ppirr}
    triangulation containing $\geq 3$ tetrahedra.  Then $G$ cannot
    contain two distinct vertices $U$ and $V$, where $U$ and $V$ are
    joined by an edge, and where $U$ and $V$ are each joined to the
    endpoints of two distinct one-ended chains.
    This prohibited situation is illustrated in Figure~\ref{fig-square}.

    \begin{figure}[htb]
        \psfrag{U}{{\small $U$}} \psfrag{V}{{\small $V$}}
        \centerline{\includegraphics[scale=0.7]{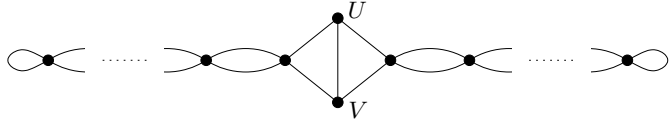}}
        \caption{Two adjacent vertices joined to two one-ended chains}
        \label{fig-square}
    \end{figure}
\end{theorem}

\begin{theorem} \label{t-mountains}
    Let $G$ be the face pairing graph of a closed minimal {\ppirr}
    triangulation containing $\geq 3$ tetrahedra.  Then $G$ cannot
    contain two distinct vertices $U$ and $V$, where $U$ and $V$ are
    each joined to the endpoints of three distinct one-ended chains.
    This prohibited situation is illustrated in Figure~\ref{fig-mountains}.

    \begin{figure}[htb]
        \psfrag{U}{{\small $U$}} \psfrag{V}{{\small $V$}}
        \centerline{\includegraphics[scale=0.7]{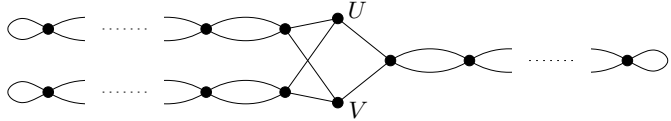}}
        \caption{Two vertices joined to three one-ended chains}
        \label{fig-mountains}
    \end{figure}
\end{theorem}

As indicated above, a certain amount of logic can be applied to both of
these theorems simultaneously.  In particular, we can obtain partial
results by analysing the graph obtained by removing the central edge in
Theorem~\ref{t-square}, or by removing the third one-ended chain in
Theorem~\ref{t-mountains}.

\begin{lemma} \label{l-emptysquare}
    Let $G$ be the face pairing graph of a closed minimal {\ppirr}
    triangulation containing $\geq 3$ tetrahedra.  Suppose that $G$
    contains two distinct vertices $U$ and $V$, where $U$ and $V$ are
    each joined to the endpoints of two distinct one-ended chains
    as illustrated in Figure~\ref{fig-emptysquare}.

    \begin{figure}[htb]
        \psfrag{U}{{\small $U$}} \psfrag{V}{{\small $V$}}
        \centerline{\includegraphics[scale=0.7]{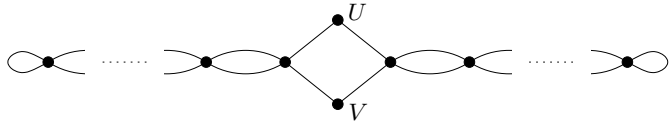}}
        \caption{Two vertices joined to two one-ended chains}
        \label{fig-emptysquare}
    \end{figure}

    Then the tetrahedra corresponding to the vertices of this subgraph
    form the structure illustrated in Figure~\ref{fig-emptysquare-tets}.
    Layered solid tori corresponding to the two chains sit beneath the
    diagram; the boundary of the first solid torus is formed from
    faces \polytri{ABE} and \polytri{BEF}, and the boundary of the
    second solid torus is formed from faces \polytri{BCF} and
    \polytri{ADE} (note the directed edge identifications
    $\dedge{AD}=\dedge{CF}$, $\dedge{AB}=\dedge{EF}$,
    $\dedge{BC}=\dedge{DE}$, and $\dedge{AE}=\dedge{BF}$).  The
    tetrahedron corresponding to graph vertex $U$ is $\mathit{ABDE}$,
    and the tetrahedron corresponding to graph vertex $V$ is
    $\mathit{BCEF}$ (note that each of these tetrahedra shares a
    face with each of the aforementioned boundary tori).

    \begin{figure}[htb]
        \psfrag{A}{{\small $A$}} \psfrag{B}{{\small $B$}}
        \psfrag{C}{{\small $C$}} \psfrag{D}{{\small $D$}}
        \psfrag{E}{{\small $E$}} \psfrag{F}{{\small $F$}}
        \psfrag{First torus boundary}{{\small First torus boundary}}
        \psfrag{Second torus boundary}{{\small Second torus boundary}}
        \centerline{\includegraphics[scale=0.6]{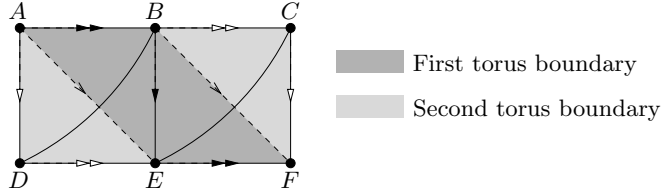}}
        \caption{The configuration of tetrahedra obtained in
            Lemma~\ref{l-emptysquare}}
        \label{fig-emptysquare-tets}
    \end{figure}
\end{lemma}

\begin{proof}
    From Theorem~\ref{t-chain-lst} we know that each one-ended chain
    corresponds to a layered solid torus in the triangulation.  Let
    $\Delta_V$ be the tetrahedron corresponding to vertex $V$ in the
    face pairing graph; since $V$ is joined to each one-ended chain by
    an edge, it follows that $\Delta_V$ is joined to each layered solid
    torus along a face.  By symmetry we may assume the leftmost diagram
    of Figure~\ref{fig-emptysquare-start}, where the first solid torus
    has boundary faces \polytri{ABE} and \polytri{BEF}, the second solid
    torus has boundary faces \polytri{BCF} and \polytri{CFG}, and
    tetrahedron $\Delta_V$ has vertices $B$, $C$, $E$ and $F$.  As in
    the theorem statement, both solid tori sit beneath the diagram.

    \begin{figure}[htb]
        \psfrag{A}{{\small $A$}} \psfrag{B}{{\small $B$}}
        \psfrag{C}{{\small $C$}} \psfrag{D}{{\small $D$}}
        \psfrag{E}{{\small $E$}} \psfrag{F}{{\small $F$}}
        \psfrag{G}{{\small $G$}} \psfrag{H}{{\small $H$}}
        \psfrag{Label1}{{\small Two solid tori plus one tetrahedron}}
        \psfrag{Label2}{{\small Two solid tori plus two tetrahedra}}
        \centerline{\includegraphics[scale=0.6]{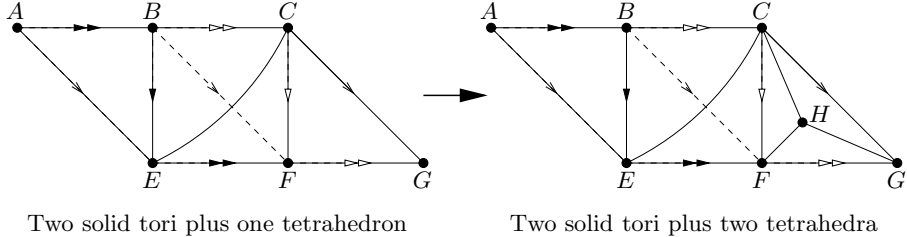}}
        \caption{Joining two tetrahedra to two layered solid tori}
        \label{fig-emptysquare-start}
    \end{figure}

    Similarly, let $\Delta_U$ be the tetrahedron corresponding to the
    graph vertex $U$.  Since $U$ is joined to each one-ended chain,
    $\Delta_U$ must be attached to both faces \polytri{ABE} and \polytri{CFG}.
    In the rightmost diagram of Figure~\ref{fig-emptysquare-start} we
    see it attached to face \polytri{CFG}, where $\Delta_U$ uses
    vertices $C$, $F$, $G$ and $H$.

    It remains to identify \polytri{ABE} with one of the three remaining
    faces of $\Delta_U$ (specifically \polytri{CHF}, \polytri{FHG} or
    \polytri{GHC}).  Most possible identifications can be eliminated by
    examining edge \polyedge{AE}, which (through its existing
    identification with \polyedge{CG})
    forms part of the torus $\mathit{BCGF}$.
    Corollary~\ref{c-torusdistinct} shows that \polyedge{AE} cannot be
    identified with \polyedge{CF} or \polyedge{FG}, and
    Corollary~\ref{c-toruslayer} shows that \polyedge{AE} cannot be
    identified with \polyedge{CH}, \polyedge{FH} or \polyedge{GH}.
    Finally, directed edge $\dedge{AE}$ cannot be identified with
    $\dedge{GC}$ since it is already identified with the reverse
    edge $\dedge{CG}$.

    The only possible identification of faces that does not involve
    one of these impossible edge matchings is for \polytri{ABE} to
    be identified with \polytri{CHG} (in particular, with directed
    edge $\dedge{AE}$ mapping to $\dedge{CG}$).  This gives us the
    configuration described in the lemma statement, and the proof is
    complete.
\end{proof}

Having achieved the partial results described in
Lemma~\ref{l-emptysquare}, we can apply these results in different
ways to obtain individual proofs of
Theorems~\ref{t-square} and~\ref{t-mountains}.

{ 
\renewcommand{\prooflabel}{Proof of Theorem~\ref{t-square}}
\begin{proof}
    Suppose the face pairing graph $G$ contains a subgraph as
    illustrated in Figure~\ref{fig-square}.  From
    Lemma~\ref{l-emptysquare} it is clear that the triangulation
    contains the structure illustrated in
    Figure~\ref{fig-emptysquare-tets}, where one additional face of
    tetrahedron $\Delta_U=\mathit{ABDE}$ is identified with one
    additional face of tetrahedron $\Delta_V=\mathit{BCEF}$.

    By symmetry we may assume the face of $\Delta_U$ is \polytri{BDE}.
    This identification of faces may cause edge \polyedge{BE} to be
    identified with some other edge in the diagram.  However, since
    \polyedge{BE} belongs to the torus $\mathit{ABEF}$, a little
    investigation shows that Corollaries~\ref{c-torusdistinct}
    and~\ref{c-toruslayer} do not allow \polyedge{BE} to be identified
    with any other edge.  It follows that the two
    identified faces must be \polytri{BDE} and \polytri{BCE}, with
    directed edge $\dedge{BE}$ mapping to itself.

    This however implies that edges \polyedge{DE} and \polyedge{CE}
    become identified, so that face \polytri{CEF} contains two distinct
    edges of the second torus (the torus with faces \polytri{BCF} and
    \polytri{ADE}).  Applying Corollary~\ref{c-toruslayer} once more
    gives us a contradiction and the proof is complete.
\end{proof}
} 

{ 
\renewcommand{\prooflabel}{Proof of Theorem~\ref{t-mountains}}
\begin{proof}
    Suppose the face pairing graph $G$ contains a subgraph as
    illustrated in Figure~\ref{fig-mountains}.  This time
    Lemma~\ref{l-emptysquare} shows that the tetrahedra form the
    structure illustrated in Figure~\ref{fig-emptysquare-tets}, where an
    additional layered solid torus $T$ is attached to one face of
    tetrahedron $\Delta_U=\mathit{ABDE}$ (call this face $F'$)
    and one face of tetrahedron $\Delta_V=\mathit{BCEF}$ (call this face $F''$).

    Attaching a layered solid torus to faces $F'$ and $F''$ will cause
    the three edges of $F'$ to be identified with the three edges of $F''$.
    By symmetry we may assume $F'=\Delta\mathit{BDE}$.  As in the
    previous proof, edge \polyedge{BE} may not be identified with any
    other edge in the diagram.  It therefore follows that
    $F''=\Delta\mathit{BCE}$, where the new solid torus gives directed
    edge identifications $\dedge{DE}=\dedge{BC}$ and $\dedge{BD}=\dedge{CE}$.
    This is illustrated in Figure~\ref{fig-mountains-threetori}, where
    the first two solid tori sit beneath the diagram as before, and the
    new solid torus sits above the diagram (but is not drawn).

    \begin{figure}[htb]
        \psfrag{A}{{\small $A$}} \psfrag{B}{{\small $B$}}
        \psfrag{C}{{\small $C$}} \psfrag{D}{{\small $D$}}
        \psfrag{E}{{\small $E$}} \psfrag{F}{{\small $F$}}
        \psfrag{First torus boundary}{{\small First torus boundary}}
        \psfrag{Second torus boundary}{{\small Second torus boundary}}
        \psfrag{Third torus boundary}{{\small Third torus boundary}}
        \centerline{\includegraphics[scale=0.6]{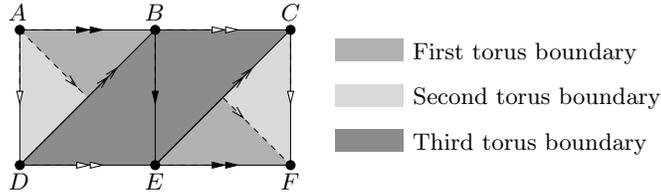}}
        \caption{The configuration of tetrahedra obtained in
            the proof of Theorem~\ref{t-mountains}}
        \label{fig-mountains-threetori}
    \end{figure}

    We conclude the proof with an examination of vertex links.  Note
    that each layered solid torus is a one-vertex structure, where the
    link of this single vertex is a disc \cite{burton-nor7,0-efficiency}.
    Moreover, our entire construction as seen in
    Figure~\ref{fig-mountains-threetori} has only one vertex, since each
    layered solid torus has only one vertex, and since vertices
    $A$, $B$, $C$, $D$, $E$ and $F$ are all identified.

    An examination of the link of this single
    vertex shows that it is not a disc ---
    instead the vertex link in our construction is a punctured torus.
    There is no way that we can complete this triangulation to give
    a 2-sphere vertex link (as required for a closed 3-manifold
    triangulation), and so the scenario described is impossible.
\end{proof}
} 

\subsection{Experimental Results} \label{s-graphs-data}

As outlined in Section~\ref{s-algm-gen}, we can use
Theorems~\ref{t-straybigon}, \ref{t-square} and~\ref{t-mountains} to
improve the generation of triangulations in the census algorithm.
Specifically, we refuse to process any face pairing graphs that do
not satisfy the constraints of these theorems.

Since the generation of face pairing graphs takes negligible time,
a very rough estimate suggests that if we can eliminate $k\%$ of graphs
in this way, we should save approximately $k\%$ of the overall running
time.  We therefore focus our efforts on measuring how many graphs are
eliminated in comparison to the total number of graphs available.

An initial set of figures is provided in
Table~\ref{tab-graphs-census} for a sense of scale.
The ``Total'' column measures the
total number of connected 4-valent graphs on the given number of
vertices.  The ``Census'' column indicates how many of these in fact
appear as face pairing graphs of closed minimal {\ppirr} triangulations
in the orientable or non-orientable census, and the ``Not used'' column
counts the graphs that do not appear (so the second column is the sum
of the final two).

\begin{table}[htb]
\caption{Total numbers of face pairing graphs on $\leq 10$ vertices}
\begin{center} \begin{tabular}{|r|r|r|r|}
    \hline
    \bf Vertices & \bf Total & \bf Census & \bf Not used \\
    \hline
    3  &       4 &   2 &       2 \\
    4  &      10 &   4 &       6 \\
    5  &      28 &   8 &      20 \\
    6  &      97 &  14 &      83 \\
    7  &     359 &  29 &     330 \\
    8  &  1\,635 &  66 &  1\,569 \\
    9  &  8\,296 & 143 &  8\,153 \\
    10 & 48\,432 & 404 & 48\,028 \\
    \hline
\end{tabular} \end{center}
\label{tab-graphs-census}
\end{table}

It should be noted that, as the number of
vertices grows large, the number of graphs that actually appear in the
census becomes negligible --- by the time we reach 10 vertices, the
number of census graphs is under $1\%$.  It is clear that a great deal
of running time could be saved if we were able to eliminate all of these
unused graphs through theorems such as those proven here.

Alas we do not achieve our dream of eliminating every unnecessary graph,
but we do manage to eliminate well over $80\%$ of them.
Table~\ref{tab-graphs-elim} shows how well our new theorems perform.
The ``Old results'' column counts the number of graphs eliminated
using the earlier Theorem~\ref{t-badgraphsold} from \cite{burton-facegraphs},
and the ``New results'' column counts the number of graphs eliminated
using the newly proven Theorems~\ref{t-straybigon}, \ref{t-square}
and~\ref{t-mountains}.  The column labelled ``All results'' combines
these figures to count the total number of graphs that can be
removed using any of the old or new results, and the ``Left over''
column counts the number of graphs that do not appear in the census but
were not eliminated at any stage.

\begin{table}[htb]
\caption{Elimination of face pairing graphs on $\leq 10$ vertices}
\begin{center} \begin{tabular}{|r|r@{~~}r|r@{~~}r|r@{~~}r|r@{~~}r|}
    \hline &
    \multicolumn{6}{|c|}{\bf Graphs eliminated} & & \\
    \cline{2-7} \bf Vertices &
    \multicolumn{2}{|c|}{\bf Old results} &
    \multicolumn{2}{|c|}{\bf New results} &
    \multicolumn{2}{|c|}{\bf All results} &
    \multicolumn{2}{|c|}{\bf Left over} \\
    \hline
    3  &       2 & (100\%) &      1 & (50\%) &       2 & (100\%) &     0 & (0\%) \\
    4  &       6 & (100\%) &      4 & (67\%) &       6 & (100\%) &     0 & (0\%) \\
    5  &      16 & (80\%) &      14 & (70\%) &      19 & (95\%) &      1 & (5\%) \\
    6  &      58 & (70\%) &      60 & (72\%) &      74 & (89\%) &      9 & (11\%) \\
    7  &     221 & (67\%) &     238 & (72\%) &     290 & (88\%) &     40 & (12\%) \\
    8  &     997 & (64\%) &  1\,116 & (71\%) &  1\,343 & (86\%) &    226 & (14\%) \\
    9  &  4\,930 & (60\%) &  5\,834 & (72\%) &  6\,904 & (85\%) & 1\,249 & (15\%) \\
    10 & 27\,681 & (58\%) & 34\,452 & (72\%) & 40\,353 & (84\%) & 7\,675 & (16\%) \\
    \hline
\end{tabular} \end{center}
\label{tab-graphs-elim}
\end{table}

Although the new results appear weak for very small numbers of vertices,
they quickly outperform the earlier results of \cite{burton-facegraphs}.
By the time we reach 10~tetrahedra they have climbed to a steady level
of eliminating 71--72\% of all unused graphs, and when combined with the
earlier results they eliminate well over $80\%$.  It should be noted
that the percentages in this table are with respect to the total
number of unused graphs (e.g., $48\,028$ graphs for 10~vertices).  However,
since the number of census graphs is negligible for larger numbers of
vertices, we should still expect to save ourselves
well over $80\%$ of the total running time.

\begin{figure}[htb]
    \centerline{\includegraphics[scale=0.6]{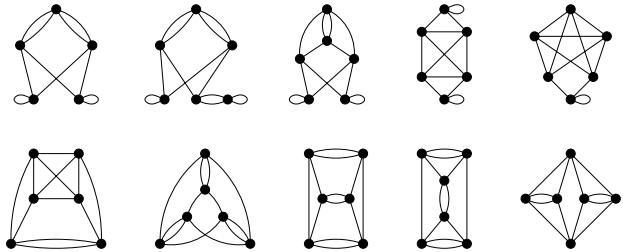}}
    \caption{Unused graphs on $\leq 6$ vertices that are not eliminated}
    \label{fig-leftover}
\end{figure}

Looking forward, the data to focus on is the final ``Left over''
column.  This lists the number of graphs that {\em could}
have been eliminated (since they do not appear in the census), but were
not eliminated by any of our old or new theorems.  Although small with
respect to the total number of graphs, these figures are not negligible
--- it is clear that with further work we could still make significant
improvements to the running time.  For reference,
Figure~\ref{fig-leftover} lists the ten unused graphs on $\leq 6$ vertices
that are not eliminated by any of theorems discussed here.

As a final note, Table~\ref{tab-graphs-freq} shows how each of the three
new theorems performs individually.\footnote{Note that the ``Total graphs''
column is generally less than the sum of the previous three, since some
graphs may be eliminated using more than one theorem.}
It is clear that
Theorem~\ref{t-straybigon} is by far the most useful of the three,
though in a census that may require months or years of processing time,
the contributions of Theorems~\ref{t-square} and~\ref{t-mountains}
should not be ignored.

\begin{table}[htb]
\caption{Frequency of prohibited subgraphs on $\leq 10$ vertices}
\begin{center} \begin{tabular}{|r|r|r|r|r|}
    \hline
    \bf Vertices & \bf Thm~\ref{t-straybigon} &
    \bf Thm~\ref{t-square} & \bf Thm~\ref{t-mountains} &
    \bf Total graphs \\
    \hline
    3  & 1 & 0 & 0 & 1 \\
    4  & 4 & 1 & 0 & 4 \\
    5  & 13 & 2 & 1 & 14 \\
    6  & 56 & 5 & 2 & 60 \\
    7  & 227 & 13 & 5 & 238 \\
    8  & 1\,083 & 46 & 14 & 1\,116 \\
    9  & 5\,730 & 170 & 47 & 5\,834 \\
    10 & 34\,059 & 746 & 176 & 34\,452 \\
    \hline
\end{tabular} \end{center}
\label{tab-graphs-freq}
\end{table}

%
%
%

%% file: ufind.tex
\section{Tracking Vertex and Edge Links} \label{s-ufind}

We now put face pairing graphs behind us and move on to our second round
of algorithmic improvements.  Here we aim to prune the search tree by
tracking partially constructed vertex and edge links, with the help of
a modified union-find algorithm.

Section~\ref{s-ufind-overview} gives an overview of how the census
algorithm is modified, including which branches of the search tree we
aim to prune, and outlines some of the implementation difficulties that
arise.
In Section~\ref{s-ufind-ufind} we introduce the classical union-find
algorithm, and in Section~\ref{s-ufind-mod} we describe how it can be
modified to address the difficulties previously raised.
Finally, experimental data is presented in Section~\ref{s-ufind-timing}
to show just how well these improvements perform in practice.

\subsection{Pruning the Search Tree} \label{s-ufind-overview}

Recall Algorithm~\ref{a-generation}, in which we split the generation of
triangulations into the enumeration and then processing of face pairing
graphs.  The processing stage consumes virtually all of the running
time, and involves a lengthy search through all possible rotations and
reflections by which tetrahedron faces might be identified according to
a particular face pairing graph.

It is natural to implement this search as a recursive (depth-first) search
with backtracking.  If we have $n$ tetrahedra then we have $2n$ pairs of
faces to join; we first choose one of the six possible rotations and
reflections for the first pair, then we choose one for the second pair
and so on.  If at any stage we find ourselves in an impossible
situation, or if the six options for a particular pair are exhausted,
we backtrack to the previous pair of faces, move on to the next rotation
or reflection for that earlier pair, and then push forwards again.  In
this way all $6^{2n}$ possible rotations and reflections are covered.

A side-effect of such a backtracking algorithm is that at each stage we
have a partially-constructed triangulation --- some tetrahedron faces
are identified in pairs, and some faces are still waiting for their
gluing instructions.  If we can identify from a partially-constructed
triangulation that the finished triangulation cannot belong to our
census, we can avoid any deeper searching and backtrack immediately.
This is what is meant by ``pruning a branch of the search tree''.  Note
that the earlier we prune, the better --- if we can identify an unwanted
triangulation after joining together $k$ pairs of faces, we remove
$6^{2n-k}$ potential triangulations from the search.

This explains the mechanics of pruning; the question remains however of
how to detect an unwanted triangulation when only some of the
tetrahedron faces have been joined together.  Fortunately several of the
results of Section~\ref{s-prelim} are useful in this regard, since they
refer only to small regions within a triangulation (e.g., single faces
or single edges).  The specific results we use include:
\begin{itemize}
    \item If an edge link has been completely constructed (i.e., the
    edge does not belong to any faces that are still not joined to their
    partners), then that edge must have degree $\geq 3$, and if it lies
    on three distinct tetrahedra then it must have degree $>3$
    (Lemma~\ref{l-edgedeg}).
    \item No vertex link may be completely constructed unless the entire
    triangulation is finished, since otherwise we will obtain a
    $\geq 2$-vertex triangulation (Lemma~\ref{l-onevertex}).
    \item No face may have two of its edges identified to form a cone,
    and no face may have all three edges identified
    (Lemmas~\ref{l-cone} and~\ref{l-l31}).
    \item No edge may be joined to itself in reverse, and no vertex may
    have a non-orientable link (required for a 3-manifold triangulation).
\end{itemize}

Whilst the following two results are global rather than local
properties, we add them to the list since they will prove easy to test
within the infrastructure that we develop.  Both results use the fact
that each pairing of tetrahedron faces introduces at most three new
edge identifications and three new vertex identifications.
\begin{itemize}
    \item If we have joined $k$ pairs of tetrahedron faces,
    the number of distinct
    vertex classes must be $\leq 1+3(2n-k)$ (Lemma~\ref{l-onevertex}).
    \item If we have joined $k$ pairs of tetrahedron faces,
    the number of distinct
    edge classes must be $\geq n+1$ and $\leq n+1+3(2n-k)$
    (Corollary~\ref{c-edgecount}).
\end{itemize}

All of the properties described above
can be tested by examining edge and vertex links
within the partially constructed triangulation, i.e., by examining which
individual tetrahedron edges or vertices come together to form a single edge
or vertex of the overall triangulation.  The simplest way of doing this
is to group the individual tetrahedron edges and vertices into
equivalence classes, where an equivalence class represents all edges or
vertices that are identified together.

We could potentially recreate these equivalence classes each time we
test our pruning constraints.  However, this would be far too slow ---
the creation of equivalence classes could take up to linear time in $n$,
and this would need to be done every time we joined a pair of
tetrahedron faces together or pulled them apart again.  A preferable solution
is to dynamically modify the equivalence classes as we modify the
triangulation, in the hope that only a small amount of work needs to be
done at each step to keep the two in sync.

The question still remains of how to implement our equivalence classes.
It is typically seen in computationally intensive problems such as this that
a careful choice of underlying data structures and algorithms can have a
great effect upon the overall running time.  This choice is the primary
focus of Sections~\ref{s-ufind-ufind} and~\ref{s-ufind-mod}.  In the
meantime, we examine the specific operations that our implementation of
equivalence classes must support.
\begin{enumerate}[(i)]
    \item {\em Merging together two equivalence classes:} \label{op-merge}
    This occurs
    when a new pair of tetrahedron faces are joined (specifically, each
    join results in $\leq 3$ merges of edge classes and $\leq 3$ merges
    of vertex classes).
    \item {\em Splitting equivalence classes:} \label{op-split}
    This occurs when we
    backtrack and pull a pair of tetrahedron faces apart.  As a result
    we must undo whatever merges were performed in the previous step.
    \item {\em Testing for equivalence:} \label{op-equiv}
    We must be able to identify
    whether two tetrahedron edges or vertices belong to the same
    equivalence class (this is required when testing for conditions such as
    faces with multiple edges identified).
    \item {\em Testing for completeness:} \label{op-complete}
    We need to identify whether an
    edge or vertex link has been completely constructed, i.e., does not
    involve any faces not yet joined to their partners (this is
    required for the one-vertex test and the edge degree tests).
    \item {\em Measuring class sizes:} \label{op-size}
    We must be able to count the
    total number of edges in an equivalence class (this is required for
    the edge degree tests).
    \item {\em Tracking orientation:} \label{op-orient}
    We must keep track
    of the orientation of tetrahedron edges within edge classes, as well
    as the orientation of the triangular pieces of vertex link that each
    tetrahedron supplies for each vertex class (this is required for
    the 3-manifold tests).
    \item {\em Counting classes:} \label{op-count}
    We must be able to count the total
    number of edge and vertex classes (this is required for the edge and
    vertex counting tests).
\end{enumerate}

It is simple enough to implement all of the above operations in a
na\"ive way.  However, all of these operations are performed
{\em extremely} frequently --- we merge and split classes every time we
step forwards and backwards in the search tree, and we run the remaining
tests every time we have a new partially constructed triangulation.
Because of this, every one of these operations must be very fast;
otherwise the time we save in pruning the search tree may be lost in
testing whether we are able to prune.

Na\"ive array-based and list-based implementations of equivalence classes
typically require linear time (or worse) for one or more of the operations
in this list.  In the following section we introduce the union-find
algorithm, and show how it can be modified to perform all of the above
operations in logarithmic time at worst.

\subsection{Introducing the Union-Find Algorithm} \label{s-ufind-ufind}

The {\em union-find} algorithm is a classical algorithm for building a
set of equivalence classes from a collection of objects and
relationships.\footnote{This algorithm frequently appears in the
equivalent setting of building a set of connected graph components from
a collection of vertices and edges.}
It is designed to perform the following operations extremely quickly:
\begin{itemize}
    \item Merging two equivalence classes, i.e.,
    operation~(\ref{op-merge}) in the list above;
    \item Testing whether two objects are equivalent, i.e.,
    operation~(\ref{op-equiv}) in the list above.
\end{itemize}

A key property of union-find is that operations of these two
types may be interspersed.  For instance, we may merge some classes,
test some pairs of objects, merge some more classes, test additional
pairs of objects, and so on.  This is particularly important for our
application, since the backtracking nature of the census algorithm means
that updates (merging and splitting classes) and tests will be
interleaved in this way.

Because of the highly optimised nature of the union-find algorithm,
it is {\em not}
well suited for other operations such as splitting classes or listing
all objects within a class --- as in many classical algorithms, the
benefit of improved speed comes at the cost of reduced functionality.
In Section~\ref{s-ufind-mod} we modify the algorithm, trading a small (but
acceptable) loss of speed for the increased functionality required by
operations (\ref{op-merge})--(\ref{op-count}) in the list above.

In the meantime however, we give a brief overview of a typical
union-find algorithm.  This algorithm is described more thoroughly
by Sedgewick \cite{sedgewick}, and Cormen et~al.~offer a detailed
discussion of its time complexity \cite{cormen-algorithms}.

\begin{algorithm}[Union-Find] \label{a-ufind}
    Suppose we have a set of objects $\{x_1,\ldots,x_n\}$ and a set of
    relationships of the form $x_i \equiv x_j$.  These relationships
    group the objects into disjoint equivalence classes.  For each
    $x_i$, denote the class containing $x_i$ by $[x_i]$.

    We store each equivalence class as a tree structure, where each
    object $x_i$ is a node in the relevant tree.  Unlike
    ordinary trees in which parents keep lists of their children, we
    only require that each child keeps track of its parent (so it is
    easy to move up the tree, but difficult to move down).
    Such a structure is illustrated in the leftmost diagram
    of Figure~\ref{fig-ufind}.

    \begin{figure}[htb]
        \psfrag{1}{{\small $x_1$}} \psfrag{2}{{\small $x_2$}}
        \psfrag{3}{{\small $x_3$}} \psfrag{4}{{\small $x_4$}}
        \psfrag{5}{{\small $x_5$}} \psfrag{6}{{\small $x_6$}}
        \psfrag{7}{{\small $x_7$}} \psfrag{8}{{\small $x_8$}}
        \psfrag{9}{{\small $x_9$}} \psfrag{10}{{\small $x_{10}$}}
        \psfrag{11}{{\small $x_{11}$}} \psfrag{12}{{\small $x_{12}$}}
        \psfrag{13}{{\small $x_{13}$}} \psfrag{14}{{\small $x_{14}$}}
        \centerline{\includegraphics[scale=0.6]{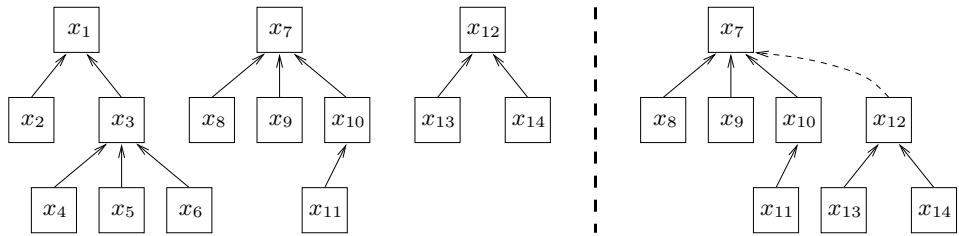}}
        \caption{Tree structures in a union-find algorithm}
        \label{fig-ufind}
    \end{figure}

    For each equivalence class, the node at the top of the tree is
    considered the {\em representative node} of that class.  For any
    object $x_i$, it is easy to locate the representative of class
    $[x_i]$ --- simply continue moving to parent nodes
    until we reach a node that has no parent.  We denote this
    representative node by $c(x_i)$.

    To test whether two objects $x_i$ and $x_j$ are equivalent, we
    simply calculate the representatives $c(x_i)$ and $c(x_j)$ as
    described above and see whether they are the same.

    To add a relationship $x_i \equiv x_j$, we again calculate the
    representatives $c(x_i)$ and $c(x_j)$.  If these representatives are
    identical then the relationship is redundant, and we discard it.
    Otherwise we merge the trees together by inserting $c(x_i)$ as a
    child of $c(x_j)$, or vice versa.  This procedure is illustrated in
    the rightmost diagram of Figure~\ref{fig-ufind}.
\end{algorithm}

It is clear that the time complexity of both adding a relationship and
testing for equivalence is bounded by the maximum tree depth.  As
it stands, the simple union-find described above can produce worst-case
trees whose depth is linear in $n$, resulting in linear time merging and
equivalence testing.  We now describe two classical
optimisations that reduce this
time to ``almost constant'', in a sense that is described more precisely
below.

The first optimisation is a form of balancing, and the second is a form
of tree compression.  It should be noted that there are several variants
of each type of optimisation; see the aforementioned references
\cite{cormen-algorithms,sedgewick} for other examples.

\begin{algorithm}[Height Balancing] \label{a-ufind-bal}
    This is a simple modification to the merging operation in a
    union-find algorithm.  For each node $x_i$, we store a number
    $d(x_i)$ which is an upper bound on the depth of the subtree rooted at
    $x_i$.  At the beginning of the algorithm (before any relationships
    are introduced), each $d(x_i)$ is set to 1.

    Suppose we are adding the relationship $x_i \equiv x_j$, which
    requires us to merge trees $c(x_i)$ and $c(x_j)$ as described in
    Algorithm~\ref{a-ufind}.  If $d(c(x_i)) < d(c(x_j))$ then we insert
    $c(x_i)$ as a child of $c(x_j)$.  Otherwise we insert
    $c(x_j)$ as a child of $c(x_i)$, and if the depths are equal
    then we increment $d(c(x_i))$ by one.
\end{algorithm}

Put simply, height balancing ensures that whenever a merge occurs, the
shorter tree is inserted as a child within the taller tree and not the
other way around.  This reduces the maximum depth to $O(\log n)$, and
therefore ensures logarithmic time for both merging and equivalence testing.

Although $d(x_i)$ is defined as an upper bound on the subtree depth, it
is clear from Algorithm~\ref{a-ufind-bal} that it is in fact the
precise subtree depth.  The reason for defining it more loosely as an
upper bound is so that height balancing can be used in conjunction with
path compression, as described below.

\begin{algorithm}[Path Compression] \label{a-ufind-comp}
    This is a modification to the union-find tree traversal operation,
    where we begin with an object $x_i$ and move upwards to locate the
    representative $c(x_i)$.  Each time we traverse the tree in this
    way, we return to $x_i$ and make a second run up the tree, this time
    changing the parent of every intermediate node to point
    directly to the root $c(x_i)$.
\end{algorithm}

This has the effect of flattening sections of trees at every convenient
opportunity, again with the aim of keeping the depths of nodes as small
as possible.  When path compression is combined with height balancing,
the running time for both merging and equivalence testing becomes
``almost constant''.
More specifically, if we perform $k$ operations on
$n$ objects then the total running time is guaranteed to be
$O(k\alpha(n))$, where $\alpha$ is a function that grows so
slowly that it can be considered constant for all practical purposes.
This result is due to Tarjan; see \cite{cormen-algorithms} for a summary
of the proof.

It should be noted that, when we compress the tree, we do not alter any
of the depth estimates $d(x_j)$.  To do so would (in the worst case)
require linear time, losing the efficiency that we have gained.
However, since compression can never increase the depth of a subtree,
it is clear that each $d(x_j)$ remains an upper bound as required by
Algorithm~\ref{a-ufind-bal}.

\subsection{Modifying the Union-Find Algorithm} \label{s-ufind-mod}

Whilst the optimised union-find algorithm is well known for its speed in
building equivalence classes, it suffers from the fact that it is
{\em only} designed for building classes --- it is not
well suited to pulling them apart again.
This makes it difficult to apply union-find to the census algorithm,
which includes backtracking as a core requirement.

Fortunately we can work around this problem with only a small reduction
in speed.  Here we describe a series of modifications that can be
applied to union-find, resulting in an algorithm that supports not
only the merging and splitting of classes, but all of the required operations
(\ref{op-merge})--(\ref{op-count}) from Section~\ref{s-ufind-overview}.

\begin{description}
    \item[Modification A (No path compression)] ~ \par
    We optimise union-find using height balancing as before,
    but path compression is abandoned.  The reason is that the
    compression of trees is a very difficult operation to undo when
    backtracking.

    By abandoning path compression, our trees become a little deeper
    on average.  However, the height
    balancing still ensures that the maximum depth is $O(\log n)$, and
    so merging and equivalence testing both remain logarithmic time
    operations.

    \item[Modification B (Support splitting)] ~ \par
    For each relationship that we add to the system,
    we store (i) whether that relationship caused two classes to be
    merged, and (ii) if so, which node formed the root of the smaller tree
    that was merged into the larger tree.  For the merge depicted in
    Figure~\ref{fig-ufind}, the node stored would be $x_{12}$.

    This makes splitting classes a constant time operation.  As we
    backtrack and remove a relationship from the system, we simply
    look up the corresponding node in our table and remove its parent
    link.  This will revert the entire set of trees to its previous state.
    To undo the merge depicted in Figure~\ref{fig-ufind}, we would look
    up $x_{12}$ in our table and remove its parent link, which is
    represented by the dotted arrow.

    Note that this procedure requires us to remove relationships in the
    reverse order to which we added them.  However, this requirement is
    trivially satisfied in the backtracking algorithm of
    Section~\ref{s-ufind-overview}.

    \item[Modification C (Add data to nodes)] ~ \par
    For each node $x_i$, we store some additional pieces of
    information to assist with our various pruning tests.
    If $x_i$ is the representative of its equivalence
    class (i.e., the root of its tree), then these pieces of information
    are as follows.
    \begin{itemize}
        \item We store the total number of objects in the equivalence class,
        denoted by $n(x_i)$;
        \item If we are examining vertex links, we also store the number of
        boundary edges in the partially-constructed link of vertex $x_i$,
        denoted by $b(x_i)$.
    \end{itemize}

    Storing $n(x_i)$ allows us to efficiently measure class sizes
    (this becomes a constant time lookup).  Likewise, storing $b(x_i)$
    allows us to efficiently test for complete vertex links, since the link
    of a representative vertex $x_i$ is complete if and only if $b(x_i)=0$.
    We do not bother storing boundary information
    in the case of edge links, since
    testing for completeness is easy --- an edge link becomes complete as soon
    as we find a new relationship between two edges already in the same
    equivalence class.

    Recall that we only defined $n(x_i)$ and $b(x_i)$ for the root node
    of each equivalence class.  We also store values of $n(x_i)$ and
    $b(x_i)$ for other child nodes, though these values differ from the roots
    and do not describe any current equivalence classes.  Instead they
    contain historical data, which we will shortly see can be used to
    support backtracking.

    When two classes are merged, only the information at the new root node
    is updated.  Specifically, suppose that trees rooted at $x_p$ and
    $x_c$ are merged, with $x_c$ inserted as a child of $x_p$.  Then
    $n(x_p)$ is replaced with $n(x_p)+n(x_c)$, since the size of the
    new tree is the sum of the two original tree sizes.  Likewise, in
    the case of vertex links, $b(x_p)$ is replaced with
    $b(x_p)+b(x_c)-2$.  Note that we subtract two because, when vertices
    $x_p$ and $x_c$ becomes identified, one boundary edge from the link of $x_p$
    becomes joined to one boundary edge from the link of $x_c$.

    It should be observed that these updates involve only constant
    time calculations based on the information stored in each original
    root node.  An example of this merging procedure can be seen
    in Figure~\ref{fig-ufind-data}, where $x_p=x_1$ and $x_c=x_3$.

    \begin{figure}[htb]
        \centerline{\includegraphics[scale=0.95]{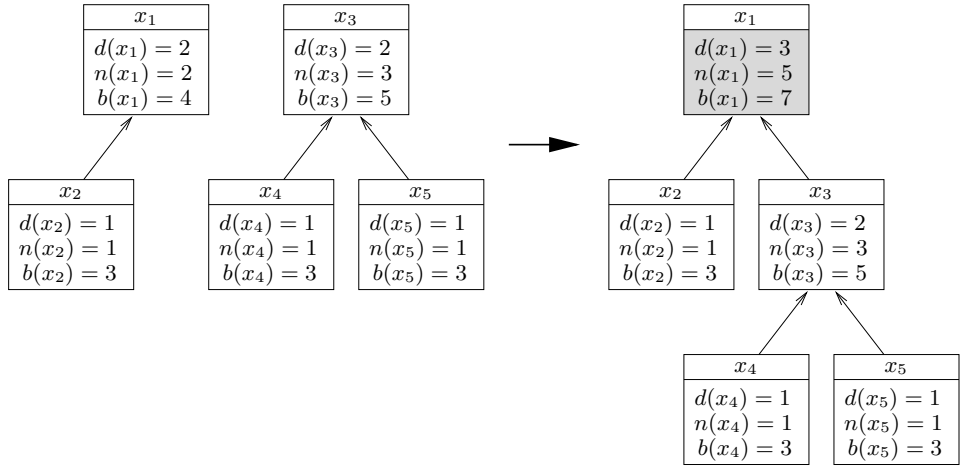}}
        \caption{Updating node data during a merge}
        \label{fig-ufind-data}
    \end{figure}

    Consider now when a class is split.  Suppose we are backtracking
    from the earlier merge operation, and are therefore cutting the
    subtree rooted at $x_c$ away from its root node $x_p$.
    Since we are restoring the system to the state just before
    the merge, we see that the data stored in
    node $x_c$ is already accurate --- it was certainly accurate before the
    merge (back when it was a root node), and it has not been changed
    since.  The old data stored in node $x_p$ is restored by replacing
    $n(x_p)$ with $n(x_p)-n(x_c)$ and replacing
    $b(x_p)$ with $b(x_p)-b(x_c)+2$.  Again these are constant time
    calculations.

    \onlyau{\newpage}

    \item[Modification D (Track orientations)] ~ \par
    For every edge of every tetrahedron in the triangulation, we
    define a canonical orientation.  Likewise, for every vertex of every
    tetrahedron, we define a canonical orientation of the corresponding
    triangular piece of vertex link.

    We use these orientations as follows.  Recall that each equivalence
    class represents a set of edges or vertices that are identified
    together.  For each child-parent relationship in a union-find tree,
    we store a boolean value (true or false)
    indicating whether the child and parent objects (edges or vertices) are
    identified in an orientation-preserving or orientation-reversing manner.

    Note that this allows us to determine in logarithmic time
    whether {\em any} two
    objects $x_i,x_j$ in the same tree are identified in an
    orientation-preserving or orientation-reversing manner --- we
    simply follow separate paths from $x_i$ and $x_j$ to the common root
    node, and keep track of these boolean values as we go.

    Consider now the situation in which two classes are merged.  This occurs
    when two objects $x_i,x_j$ from different equivalence classes become
    identified.  The merge requires us to add a new child-parent
    relationship between the root nodes $c(x_i)$ and $c(x_j)$; we must
    therefore establish whether or not this new relationship preserves
    orientation.  This is simple to determine --- we already know how the
    orientations of $x_i$ and $x_j$ relate (since this is the identification
    that is causing the merge), and we can see how $x_i$ and $c(x_i)$
    relate and how $x_j$ and $c(x_j)$ relate by following paths to the
    root of each tree as described above.

    The value of storing these booleans is that we can now
    efficiently test for non-orientable vertex links, as well as for
    tetrahedron edges that are identified with themselves in reverse
    (recall that either condition allows us to prune our search tree).
    Such conditions arise when we attempt to identify two tetrahedron
    edges or two tetrahedron vertices that already belong to the same
    equivalence class.  By following paths in the corresponding
    union-find tree from the identified objects $x_i,x_j$ to the common root
    node $c(x_i)=c(x_j)$ and tracking these boolean values,
    we can establish whether $x_i$ and $x_j$ are already
    identified in an orientation-preserving or
    orientation-reversing manner, and whether the new identification
    contradicts this.

    \item[Modification E (Count classes)] ~ \par
    Finally, we keep a separate variable that counts the total number
    of equivalence classes.  This variable is decremented each time we
    merge two classes, and incremented each time we split a class.  This
    allows us to query the total number of equivalence classes in constant
    time.
\end{description}

We can now summarise the efficiency of each of the required operations
(\ref{op-merge})--(\ref{op-count}) from Section~\ref{s-ufind-overview}.
Merging classes and locating root nodes are logarithmic time operations,
since height balancing ensures that trees have a logarithmic worst-case
depth.  Likewise, testing for equivalence is logarithmic time.
Both splitting classes and counting classes are constant time operations, and
once the class representative (i.e., the root node) has been located,
testing for completeness and measuring class sizes become constant time
also.  Finally, tracking orientations requires logarithmic time since the
relevant procedures each involve two child-to-root traversals through a
tree.

We see then that the time complexity of every required operation is
logarithmic at worst.  This is a small price to pay for such enhanced
functionality, and in the following section we see that it is definitely
worth the cost.

\subsection{Experimental Results} \label{s-ufind-timing}

We conclude this section with experimental data that measures the
running time of the census algorithm with and without the help of the
modified union-find described in Section~\ref{s-ufind-mod}.

Table~\ref{tab-ufind-timing} lists precise running times for the
generation of non-orientable census triangulations for 5--7
tetrahedra.  The range 5--7 was chosen because smaller censuses all
give running times of $\leq 5$ seconds, and because for larger censuses
it is infeasible to measure the running time of the original algorithm.

\begin{table}[htb]
\caption{Running times for the modified union-find algorithm
    ({\em hh}:{\em mm}:{\em ss})}
\begin{center} \begin{tabular}{|r|r|r|r|r|}
    \hline
    & & \multicolumn{3}{|c|}{\bf Modified algorithm} \\
    \cline{3-5}
    \bf Tetrahedra & \bf Old algorithm &
        \bf Vertex links & \bf Edge links & \bf Both \\
    \hline
    5  &         8:38 &        1:24 &    0:06 & 0:02 \\
    6  &     13:49:27 &     1:23:15 &    2:40 & 0:46 \\
    7  & $\sim 65$ days & $\sim 4$ days & 2:06:05 & 21:38 \\
    \hline
\end{tabular} \end{center}
\label{tab-ufind-timing}
\end{table}

The ``Old algorithm'' column gives the running time of the original
census algorithm.  Note that some pruning is still performed in this
case; in particular, the original algorithm performs tests that do not require
the construction of large edge links or vertex links, such as tests for
edges of degree $\leq 3$.  The remaining
three columns list the running times for the new census algorithm, according
to whether we use a modified union-find structure to track vertex links,
edge links or both.  All times are measured on a single 1.7~GHz IBM Power5
processor.

The results are spectacular.  The 5-tetrahedron census is reduced from
$8\frac12$~minutes to 2~seconds, and the 7-tetrahedron census is reduced
from 65~days to a mere 22~minutes.  This shows improvements of several
orders of magnitude, increasing the speed by over 4000 times in the
7-tetrahedron case.

It is only through the use of this modified algorithm that a
10-tetrahedron non-orientable census was even feasible.  The results of
this census (as well as the sister 10-tetrahedron orientable census) are
discussed in the following section.

%% file: census.tex
\section{Census Results} \label{s-census}

Using the algorithmic improvements of Sections~\ref{s-graphs} and~\ref{s-ufind},
it has been possible to generate new census results in both the
non-orientable and orientable cases for up to 10 tetrahedra.  In particular:
\begin{itemize}
    \item The closed non-orientable census has been extended well beyond
    its previous limit of 8~tetrahedra, with a much richer
    variety of structures appearing as a result;
    \item The closed orientable census
    has been extended from just a list of manifolds to a full list
    of all minimal triangulations of these manifolds.
\end{itemize}
The non-orientable and orientable cases are treated separately in
Sections~\ref{s-census-nor} and~\ref{s-census-or} below.

The full census data, including all $16\,109$
closed minimal {\ppirr} triangulations
formed from $\leq 10$ tetrahedra (both orientable and non-orientable), is
far too large to include here.  This data may be downloaded from the
{\regina} website \cite{regina} in the form of {\tt .rga} data files,
which may be loaded into {\regina} and examined, analysed
and/or exported to some other format.\footnote{Note that {\regina}
provides a scripting facility, which may assist with analysing
and/or exporting large numbers of triangulations in bulk.}
As of version~4.3 these data files
are shipped as part of {\regina}, and may be accessed through the
{\em File}\,$\rightarrow$\,{\em Open Example} menu entry.  Each data file
begins with a text description of the census parameters and notation used.

The physical hardware requirements for these census runs were filled
by the excellent resources of the Victorian Partnership for Advanced
Computing (VPAC).  {\regina}'s census generation routines provide
native support for parallel processing, allowing for an effective use of
VPAC's cluster environments.

\subsection{Non-Orientable Census} \label{s-census-nor}

The 7-tetrahedron non-orientable censuses of Amendola and Martelli
\cite{italian-nor7} and Burton \cite{burton-nor7} show that there are
eight distinct closed non-orientable {\ppirr} 3-manifolds that can be
formed from $\leq 7$ tetrahedra.  All of these 3-manifolds fall into one
of the following categories:
\begin{itemize}
    \item Torus bundles over the circle;
    \item {\sfslong}s over either $\rpp$ or $\discref$ with two exceptional
    fibres, where $\discref$ is the disc with reflector boundary.
\end{itemize}

Subsequent results of Burton \cite{burton-nor8} show that when the limit
is raised to $\leq 8$ tetrahedra, a total of 18 distinct 3-manifolds are
found, all of which belong to the same two categories listed above.
We see then that the 8-tetrahedron non-orientable census remains
extremely limited in scope.

Happily a much richer variety of 3-manifolds appears in the
9-tetra\-hedron and 10-tetra\-hedron non-orientable census tables.
We begin with a formal presentation of the new non-orientable census
results, and follow this with a summary of the new types
of 3-manifolds that are found.

\begin{theorem} \label{t-nor-census}
    A total of 136 distinct closed non-orientable {\ppirr} 3-manifolds
    can be constructed using $\leq 10$ tetrahedra.  These manifolds are
    listed individually in Table~\ref{tab-mfdlist} in the appendix, and are
    described by a total of 1390 distinct minimal triangulations.

    A breakdown of these figures can be found in Table~\ref{tab-summary-nor},
    which lists the number of distinct 3-manifolds and distinct
    minimal triangulations for 6, 7, 8, 9 and 10 tetrahedra.
    None of these 3-manifolds can be constructed from five tetrahedra or fewer.
\end{theorem}

\begin{table}[htb]
\caption{Closed non-orientable census results}
\begin{center} \begin{tabular}{|r|r|r|}
    \hline
    \bf Tetrahedra & \bf 3-Manifolds & \bf Triangulations \\
    \hline
    $\leq 5$ & 0 & 0 \\
    6  & 5  &   24 \\
    7  & 3  &   17 \\
    8  & 10 &   59 \\
    9  & 33 &  307 \\
    10 & 85 &  983 \\
    \hline
    \bf Total & 136 & 1\,390 \\
    \hline
\end{tabular} \end{center}
\label{tab-summary-nor}
\end{table}

It is immediate from Table~\ref{tab-mfdlist} that the 9-tetrahedron and
10-tetrahedron manifolds are greater not just in number but also in
variety.  In addition to the two original categories described above,
the following new types of 3-manifolds are also seen:
\begin{itemize}
    \item {\sfslong}s over $\rpp$ and $\discref$ with three exceptional
    fibres;
    \item {\sfslong}s over the torus, Klein bottle, $\mobref$ and
    $\annrefref$ with one exceptional fibre,
    where $\mobref$ is the {\mobius} band with reflector boundary and
    $\annrefref$ is the annulus with two reflector boundaries;
    \item Non-geometric graph manifolds of the following types:
    \begin{itemize}
        \item {\sfslong}s over the disc with two exceptional fibres
        joined to {\sfslong}s over the {\mobius} band or $\annref$
        with one exceptional fibre, where $\annref$ is the
        annulus with one reflector boundary;
        \item {\sfslong}s over the annulus whose two boundary tori are
        joined together;
        \item The special manifold $\specialnine$, described below;
    \end{itemize}
    \item The special non-geometric manifold $\specialten$, also
    described below.
\end{itemize}

Two ``special manifolds'' are referred to in the list above.  These have
noteworthy properties that set them apart from the remaining manifolds
in the census.  Each special manifold is described in turn.
\begin{itemize}
    \item {\em The 9-tetrahedron manifold $\specialnine$}:

    This is a non-geometric graph manifold, formed from a single
    {\sfslong} by identifying its two boundaries.  What sets this apart
    from the other graph manifolds is that the boundaries are Klein
    bottles, not tori.

    The base orbifold of the {\sfslong} is $\annhalfref$, the annulus
    with one regular boundary curve and one half-reflector, half-regular
    boundary curve.  In addition, the path around the annulus is a
    fibre-reversing curve in the {\sfslong} (this is the meaning of
    the symbol $o_2$).  This base orbifold is illustrated in
    the leftmost diagram of Figure~\ref{fig-special9}.

    \begin{figure}[htb]
    \psfrag{a}{{\small $a$}} \psfrag{b}{{\small $b$}} \psfrag{c}{{\small $c$}}
    \psfrag{x}{{\small $x$}} \psfrag{y}{{\small $y$}}
    \psfrag{a'}{{\small $a'$}} \psfrag{b'}{{\small $b'$}}
    \psfrag{Base annulus}{{\small Base annulus}}
    \psfrag{Outer boundary}{{\small Outer boundary}}
    \psfrag{Inner boundary}{{\small Inner boundary}}
    \psfrag{(half-ref)}{{\small (reflector)}}
    \psfrag{=}{{\small $\equiv$}}
    \centerline{\includegraphics[scale=0.7]{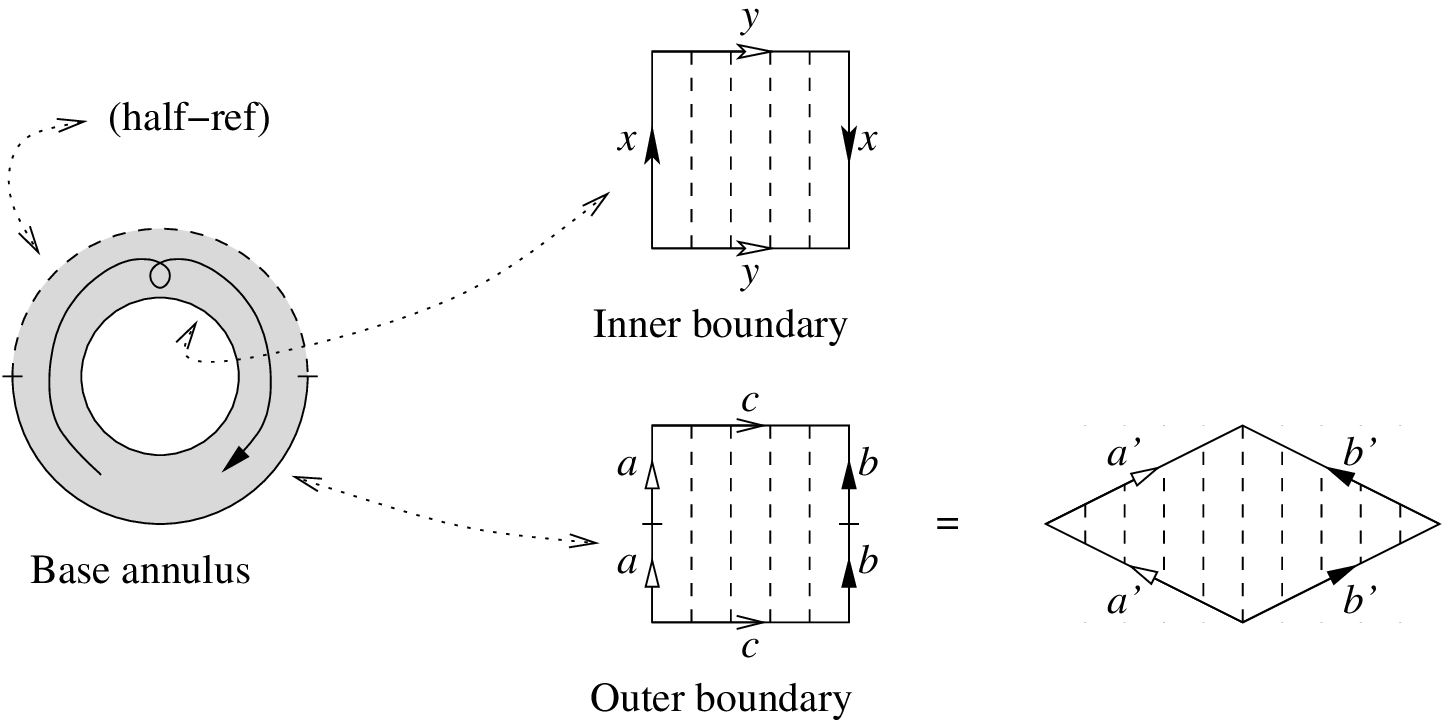}}
    \caption{The base orbifold and boundary Klein bottles for
        $\mathrm{SFS}\left(\smash{\annhalfref}/o_2\right)$}
    \label{fig-special9}
    \end{figure}

    The resulting {\sfslong} has two Klein bottle boundaries.  The
    lower central diagram of Figure~\ref{fig-special9} shows the Klein bottle
    corresponding to the half-regular boundary on the outside of the
    annulus, and the upper central diagram shows the Klein bottle
    corresponding to the regular boundary on the inside of the annulus.
    The dashed lines on each Klein bottle depict the fibres of the
    {\sfslong}.

    It is clear that the two Klein bottles cannot be identified so that
    the two sets of fibres match.  Instead, consider the diamond
    representation of the outer boundary Klein bottle as depicted
    in the rightmost diagram of Figure~\ref{fig-special9}.
    The two Klein bottles are identified so that the fibre $x$ maps to
    the transverse curve $a'b'^{-1}$, and the closed curve $a'$ (which
    is isotopic to the closed half-fibre $a$) maps to the transverse curve $y$.

    \item {\em The 10-tetrahedron manifold $\specialten$}:

    This is also non-geometric, but it is not a graph manifold.
    It consists of a Seifert fibred block and a hyperbolic block joined
    along a common Klein bottle boundary.

    The base orbifold of the Seifert fibred block is $\dischalfref$,
    the disc with half-reflector, half-regular boundary.  A single
    exceptional fibre of index~2 is also present within this block.

    The hyperbolic block is a truncated Gieseking manifold.
    The {\em Gieseking manifold} is the smallest volume cusped
    hyperbolic manifold \cite{adams-gieseking}.
    It is non-orientable with a single Klein
    bottle cusp, and has hyperbolic volume $1.01494161$.  In the
    census shipped with {\snappea} \cite{snappea} it is listed as the manifold
    {\tt m000}.

    Here we truncate the Gieseking manifold at its cusp, giving rise
    to a real Klein bottle boundary.  This boundary is then identified
    with the boundary of the Seifert fibred block, producing the closed
    manifold identified here.
\end{itemize}

Note that, with the single exception of the special 10-tetrahedron
manifold described above, every manifold in this census is a graph
manifold.  Theorem~\ref{t-nor-census} therefore yields the following
corollary.

\begin{corollary}
    All closed non-orientable {\ppirr} 3-manifolds formed from
    $\leq 9$ tetrahedra are graph manifolds.  Moreover, this bound is
    sharp --- there exists a closed non-orientable {\ppirr} 3-manifold
    formed from 10 tetrahedra whose JSJ decomposition contains a
    hyperbolic block.
\end{corollary}

It can also be noted that none of the manifolds in this census are
closed hyperbolic manifolds.  This contrasts with the orientable census,
in which closed hyperbolic manifolds first appear at the 9-tetrahedron
level \cite{matveev-complexity}.

We can nevertheless deduce the point at which closed non-orientable
hyperbolic manifolds will appear.  Hodgson and Weeks present a list of
candidate smallest volume non-orientable hyperbolic manifolds
\cite{closedhypcensus}.  The smallest volume manifold from this list
has volume $2.02988321$, and can be triangulated using 11 tetrahedra.
We therefore obtain the following sharp bound.

\begin{corollary}
    For any closed non-orientable {\ppirr} hyperbolic 3-manifold $M$,
    a minimal triangulation of $M$ must contain $\geq 11$ tetrahedra.
    Moreover, this bound is sharp --- the non-orientable
    manifold of volume $2.02988321$ from the Hodgson-Weeks census
    \cite{closedhypcensus} can be triangulated using precisely 11
    tetrahedra.
\end{corollary}

We close these non-orientable census results with a discussion of the
combinatorial structures of minimal triangulations.  In the 8-tetrahedron
non-orientable census paper \cite{burton-nor8}, the following conjectures
are made:
\begin{itemize}
    \item Every minimal triangulation of a non-flat torus bundle over
    the circle is a {\em layered torus bundle} (Conjecture~3.1);
    \item Every minimal triangulation of a non-flat {\sfslong} over
    $\rpp$ or $\discref$ with two exceptional fibres is either a
    {\em plugged thin $I$-bundle} or a {\em plugged thick $I$-bundle}
    (Conjecture~3.2).
\end{itemize}
Here layered torus bundles and plugged thin/thick $I$-bundles are
specific families of triangulations; see either \cite{burton-nor8} for an
overview or \cite{burton-nor7} for full details.

The new census results confirm these conjectures at the 9-tetrahedron and
10-tetrahedron levels.  Specifically, for 9 and 10 tetrahedra, every
minimal triangulation of a torus bundle over the circle is a layered
torus bundle, and every minimal triangulation of a {\sfslong} over
$\rpp$ or $\discref$ with two exceptional fibres is a plugged thin or
thick $I$-bundle.

\subsection{Orientable Census} \label{s-census-or}

As discussed in the introduction, closed orientable census results are
presently known for $\leq 11$ tetrahedra (though Matveev has recently
announced partial results for 12 tetrahedra also \cite{matveev12}).
However, for 7~tetrahedra and above these results are presented by
their respective authors as lists of manifolds.  Full lists of all
minimal triangulations (or even the numbers of such triangulations)
are not included.

The most comprehensive lists of orientable minimal triangulations to
date are those of Matveev \cite{matveev6}, enumerated for $\leq 6$
tetrahedra and presented in the equivalent language of special spines.
Here we extend this list of minimal triangulations to $\leq 10$
tetrahedra.

\begin{theorem}
    There are precisely $14\,719$ minimal triangulations of closed
    orientable irreducible 3-manifolds that are constructed using
    $\leq 10$ tetrahedra.  A breakdown of this count into
    $1,2,\ldots,10$ tetrahedra is provided in Table~\ref{tab-summary-or}.
\end{theorem}

\begin{table}[htb]
\caption{Closed orientable census results}
\begin{center} \begin{tabular}{|r|r|r|}
    \hline
    \bf Tetrahedra & \bf 3-Manifolds & \bf Triangulations \\
    \hline
    1  & 3 & 4 \\
    2  & 6 & 9 \\
    3  & 7 & 7 \\
    4  & 14 & 15 \\
    5  & 31 & 40 \\
    6  & 74 & 115 \\
    7  & 175 & 309 \\
    8  & 436 & 945 \\
    9  & 1\,154 & 3\,031 \\
    10 & 3\,078 & 10\,244 \\
    \hline
    \bf Total & 4\,978 & 14\,719 \\
    \hline
\end{tabular} \end{center}
\label{tab-summary-or}
\end{table}

Unlike the non-orientable census, both the triangulations and the
3-manifolds are too numerous to list here.  As described in the
beginning of Section~\ref{s-census}, the full data can be downloaded
from the {\regina} website \cite{regina}.

Note that Table~\ref{tab-summary-or} also lists the number of
distinct 3-manifolds obtained at each level of the census.
These manifold counts agree with the previous results of
Martelli \cite{italian10} and Matveev \cite{matveev11}.

Moreover, the full list of $4\,978$ manifolds obtained here with {\regina}
has been compared with the lists of Martelli and Petronio \cite{italianlists}.
A detailed matching confirms that the $4\,978$ distinct manifolds obtained in
this orientable census are identical to the $4\,978$ distinct manifolds
obtained by Martelli and Petronio.  Since both sets of census results
rely heavily on computer searches, an independent verification such as
this provides an extra level of confidence in both sets of software.

%% file: manifolds.tex
\newpage

\section*{Appendix: Tables of Non-Orientable Manifolds}

In this appendix we present the full list of all closed non-orientable
{\ppirr} 3-manifolds that can be constructed using 10 tetrahedra or
fewer, as discussed in Section~\ref{s-census-nor}.
This list, given as Table~\ref{tab-mfdlist} on the following pages,
is divided into sections
according to the number of tetrahedra in a minimal triangulation.
Note that none of these 3-manifolds can be formed from
fewer than six tetrahedra.

For each manifold $M$, this list also presents the first homology group
$H_1(M)$ and the number of distinct minimal triangulations of $M$.

\begin{table}[htb]
\caption{Notation for base orbifolds of {\sfslong}s}
\label{tab-basenotation}
\begin{center}
\begin{tabular}{l|l}
    \bf Symbol & \bf Base orbifold \\
    \hline
    $D$ & Disc \\
    $\discref$ & Disc with reflector boundary \\
    $\dischalfref$ & Disc with half-reflector, half-regular boundary \\
    $\mb$ & {\mobius} band \\
    $\mobref$ & {\mobius} band with reflector boundary \\
    $\mobref/n_2$ & {\mobius} band with reflector boundary and \\
                  & \qquad fibre-reversing curves \\
    $A$ & Annulus \\
    $\annref$ & Annulus with one reflector boundary \\
    $\annrefref$ & Annulus with two reflector boundaries \\
    $\annrefref/o_2$ & Annulus with two reflector boundaries and \\
                  & \qquad fibre-reversing curves \\
    $\annhalfref/o_2$ & Annulus with one half-reflector, half-regular \\
                      & \qquad boundary and fibre-reversing curves \\
    $\rpp$ & Projective plane \\
    $\torus/o_2$ & Torus containing fibre-reversing curves \\
    $\kb$ & Klein bottle \\
    $\kb/n_3$ & Klein bottle containing fibre-reversing curves, \\
              & \qquad but with total space remaining non-orientable \\
\end{tabular}
\end{center}
\end{table}

The notation used for the base orbifolds of {\sfslong}s is explained
in Table~\ref{tab-basenotation}.  In addition,
non-geometric graph manifolds appear in the list as follows.
\begin{itemize}
    \item $\mathrm{SFS}\left(\ldots\right) \bigcup_N
        \mathrm{SFS}\left(\ldots\right)$:

    This represents two {\sfslong}s, each with a single torus boundary,
    where these spaces are joined together along their torus boundaries
    according to the $2 \times 2$ matching matrix $N$.

    \item $\mathrm{SFS}\left(\ldots\right) / N$:

    This represents a single {\sfslong} with two torus boundaries,
    where these boundaries are identified according to the
    $2 \times 2$ matching matrix $N$.
\end{itemize}

The special non-geometric manifolds $\specialnine$ and $\specialten$ do
not fit neatly into this notation scheme; see Section~\ref{s-census-nor}
for a detailed description of each.

\newpage

\renewcommand{\arraystretch}{1.3}

\begin{center}
\tablecaption{Details for each closed non-orientable {\ppirr} 3-manifold}
\label{tab-mfdlist}
\tablehead{
    \hline
    \bf \# Tet. & \bf 3-Manifold & \bf \# Tri. & \bf Homology \\
    \hline
}
\tabletail{\hline}
\begin{supertabular}{|c|l|c|l|l|}
6 & $\torus \times I / \homtwotable{1}{1}{1}{0}$ & 1 &
    $\Z$ \\
  & $\torus \times I / \homtwotable{0}{1}{1}{0}$ & 6 &
    $\Z \oplus \Z$ \\
  & $\torus \times I / \homtwotable{1}{0}{0}{-1}$ & 3 &
    $\Z \oplus \Z \oplus \Z_2$ \\
  & $\sfs{\rpp}{(2,1)\ (2,1)}$ & 9 &
    $\Z \oplus \Z_4$ \\
  & $\sfs{\discref}{(2,1)\ (2,1)}$ & 5 &
    $\Z \oplus \Z_2 \oplus \Z_2$ \\
\hline
7 & $\torus \times I / \homtwotable{2}{1}{1}{0}$ & 4 &
    $\Z \oplus \Z_2$ \\
  & $\sfs{\rpp}{(2,1)\ (3,1)}$ & 10 &
    $\Z$ \\
  & $\sfs{\discref}{(2,1)\ (3,1)}$ & 3 &
    $\Z \oplus \Z_2$ \\
\hline
8 & $\torus \times I / \homtwotable{3}{1}{1}{0}$ & 10 &
    $\Z \oplus \Z_3$ \\
  & $\torus \times I / \homtwotable{3}{2}{2}{1}$ & 2 &
    $\Z \oplus \Z_2 \oplus \Z_2$ \\
  & $\sfs{\rpp}{(2,1)\ (4,1)}$ & 10 &
    $\Z \oplus \Z_2$ \\
  & $\sfs{\rpp}{(2,1)\ (5,2)}$ & 10 &
    $\Z$ \\
  & $\sfs{\rpp}{(3,1)\ (3,1)}$ & 7 &
    $\Z \oplus \Z_6$ \\
  & $\sfs{\rpp}{(3,1)\ (3,2)}$ & 9 &
    $\Z \oplus \Z_3$ \\
  & $\sfs{\discref}{(2,1)\ (4,1)}$ & 3 &
    $\Z \oplus \Z_2 \oplus \Z_2$ \\
  & $\sfs{\discref}{(2,1)\ (5,2)}$ & 3 &
    $\Z \oplus \Z_2$ \\
  & $\sfs{\discref}{(3,1)\ (3,1)}$ & 3 &
    $\Z \oplus \Z_3$ \\
  & $\sfs{\discref}{(3,1)\ (3,2)}$ & 2 &
    $\Z \oplus \Z_3$ \\
\hline
9 & $\torus \times I / \homtwotable{4}{1}{1}{0}$ & 18 &
    $\Z \oplus \Z_4$ \\
  & $\torus \times I / \homtwotable{4}{3}{3}{2}$ & 10 &
    $\Z \oplus \Z_6$ \\
  & $\sfs{\rpp}{(2,1)\ (5,1)}$ & 10 &
    $\Z$ \\
  & $\sfs{\rpp}{(2,1)\ (7,2)}$ & 10 &
    $\Z$ \\
  & $\sfs{\rpp}{(2,1)\ (7,3)}$ & 10 &
    $\Z$ \\
  & $\sfs{\rpp}{(2,1)\ (8,3)}$ & 10 &
    $\Z \oplus \Z_2$ \\
  & $\sfs{\rpp}{(3,1)\ (4,1)}$ & 10 &
    $\Z$ \\
  & $\sfs{\rpp}{(3,1)\ (4,3)}$ & 10 &
    $\Z$ \\
  & $\sfs{\rpp}{(3,1)\ (5,2)}$ & 10 &
    $\Z$ \\
  & $\sfs{\rpp}{(3,1)\ (5,3)}$ & 10 &
    $\Z \oplus \Z_2$ \\
  & $\sfs{\rpp}{(2,1)\ (2,1)\ (2,1)}$ & 18 &
    $\Z \oplus \Z_2 \oplus \Z_2$ \\
  & $\sfs{\discref}{(2,1)\ (5,1)}$ & 3 &
    $\Z \oplus \Z_2$ \\
  & $\sfs{\discref}{(2,1)\ (7,2)}$ & 3 &
    $\Z \oplus \Z_2$ \\
  & $\sfs{\discref}{(2,1)\ (7,3)}$ & 3 &
    $\Z \oplus \Z_2$ \\
  & $\sfs{\discref}{(2,1)\ (8,3)}$ & 3 &
    $\Z \oplus \Z_2 \oplus \Z_2$ \\
  & $\sfs{\discref}{(3,1)\ (4,1)}$ & 3 &
    $\Z \oplus \Z_2$ \\
  & $\sfs{\discref}{(3,1)\ (4,3)}$ & 3 &
    $\Z \oplus \Z_2$ \\
  & $\sfs{\discref}{(3,1)\ (5,2)}$ & 3 &
    $\Z$ \\
  & $\sfs{\discref}{(3,1)\ (5,3)}$ & 3 &
    $\Z$ \\
  & $\sfs{\discref}{(2,1)\ (2,1)\ (2,1)}$ & 5 &
    $\Z \oplus \Z_2 \oplus \Z_2 \oplus \Z_2$ \\
  & $\sfs{\torus/o_2}{(2,1)}$ & 12 &
    $\Z \oplus \Z$ \\
  & $\sfs{\kb}{(2,1)}$ & 21 &
    $\Z \oplus \Z$ \\
9 (ctd.)
  & $\sfs{\kb/n_3}{(2,1)}$ & 39 &
    $\Z \oplus \Z_8$ \\
  & $\sfs{\mobref}{(2,1)}$ & 9 &
    $\Z \oplus \Z \oplus \Z_2$ \\
  & $\sfs{\mobref/n_2}{(2,1)}$ & 9 &
    $\Z \oplus \Z_2 \oplus \Z_4$ \\
  & $\sfs{\annrefref}{(2,1)}$ & 2 &
    $\Z \oplus \Z \oplus \Z_2 \oplus \Z_2$ \\
  & $\sfs{\annrefref/o_2}{(2,1)}$ & 2 &
    $\Z \oplus \Z_2 \oplus \Z_2 \oplus \Z_2$ \\
  & $\sfs{D}{(2,1) (2,1)} \bigcup_{\homtwotable{0}{1}{1}{0}} \sfs{\mb}{(2,1)}$ & 10 &
    $\Z \oplus \Z_4$ \\
  & $\sfs{D}{(2,1) (2,1)} \bigcup_{\homtwotable{0}{1}{1}{0}} \sfs{\annref}{(2,1)}$ & 3 &
    $\Z \oplus \Z_2 \oplus \Z_4$ \\
  & $\sfs{A}{(2,1)} / {\homtwotable{0}{-1}{1}{0}}$ & 35 &
    $\Z$ \\
  & $\sfs{A}{(2,1)} / {\homtwotable{0}{-1}{1}{1}}$ & 7 &
    $\Z \oplus \Z_3$ \\
  & $\sfs{A}{(2,1)} / {\homtwotable{1}{2}{0}{1}}$ & 1 &
    $\Z \oplus \Z_2$ \\
  & $\specialnine$ & 2 &
    $\Z \oplus \Z_4$ \\
\hline
10& $\torus \times I / \homtwotable{5}{1}{1}{0}$ & 30 &
    $\Z \oplus \Z_5$ \\
  & $\torus \times I / \homtwotable{5}{4}{4}{3}$ & 16 &
    $\Z \oplus \Z_2 \oplus \Z_4$ \\
  & $\torus \times I / \homtwotable{8}{3}{3}{1}$ & 14 &
    $\Z \oplus \Z_9$ \\
  & $\torus \times I / \homtwotable{8}{5}{5}{3}$ & 2 &
    $\Z \oplus \Z_{11}$ \\
  & $\sfs{\rpp}{(2,1)\ (6,1)}$ & 10 &
    $\Z \oplus \Z_4$ \\
  & $\sfs{\rpp}{(2,1)\ (9,2)}$ & 10 &
    $\Z$ \\
  & $\sfs{\rpp}{(2,1)\ (9,4)}$ & 10 &
    $\Z$ \\
  & $\sfs{\rpp}{(2,1)\ (10,3)}$ & 10 &
    $\Z \oplus \Z_4$ \\
  & $\sfs{\rpp}{(2,1)\ (11,3)}$ & 10 &
    $\Z$ \\
  & $\sfs{\rpp}{(2,1)\ (11,4)}$ & 10 &
    $\Z$ \\
  & $\sfs{\rpp}{(2,1)\ (12,5)}$ & 10 &
    $\Z \oplus \Z_2$ \\
  & $\sfs{\rpp}{(2,1)\ (13,5)}$ & 10 &
    $\Z$ \\
  & $\sfs{\rpp}{(3,1)\ (5,1)}$ & 10 &
    $\Z \oplus \Z_2$ \\
  & $\sfs{\rpp}{(3,1)\ (5,4)}$ & 10 &
    $\Z$ \\
  & $\sfs{\rpp}{(3,1)\ (7,2)}$ & 10 &
    $\Z$ \\
  & $\sfs{\rpp}{(3,1)\ (7,3)}$ & 10 &
    $\Z \oplus \Z_2$ \\
  & $\sfs{\rpp}{(3,1)\ (7,4)}$ & 10 &
    $\Z$ \\
  & $\sfs{\rpp}{(3,1)\ (7,5)}$ & 10 &
    $\Z \oplus \Z_2$ \\
  & $\sfs{\rpp}{(3,1)\ (8,3)}$ & 10 &
    $\Z$ \\
  & $\sfs{\rpp}{(3,1)\ (8,5)}$ & 10 &
    $\Z$ \\
  & $\sfs{\rpp}{(4,1)\ (4,1)}$ & 7 &
    $\Z \oplus \Z_8$ \\
  & $\sfs{\rpp}{(4,1)\ (4,3)}$ & 9 &
    $\Z \oplus \Z_8$ \\
  & $\sfs{\rpp}{(4,1)\ (5,2)}$ & 10 &
    $\Z$ \\
  & $\sfs{\rpp}{(4,1)\ (5,3)}$ & 10 &
    $\Z$ \\
  & $\sfs{\rpp}{(5,2)\ (5,2)}$ & 7 &
    $\Z \oplus \Z_{10}$ \\
  & $\sfs{\rpp}{(5,2)\ (5,3)}$ & 9 &
    $\Z \oplus \Z_5$ \\
  & $\sfs{\rpp}{(2,1)\ (2,1)\ (3,1)}$ & 92 &
    $\Z \oplus \Z_4$ \\
  & $\sfs{\discref}{(2,1)\ (6,1)}$ & 3 &
    $\Z \oplus \Z_2 \oplus \Z_2$ \\
  & $\sfs{\discref}{(2,1)\ (9,2)}$ & 3 &
    $\Z \oplus \Z_2$ \\
10 (ctd.)
  & $\sfs{\discref}{(2,1)\ (9,4)}$ & 3 &
    $\Z \oplus \Z_2$ \\
  & $\sfs{\discref}{(2,1)\ (10,3)}$ & 3 &
    $\Z \oplus \Z_2 \oplus \Z_2$ \\
  & $\sfs{\discref}{(2,1)\ (11,3)}$ & 3 &
    $\Z \oplus \Z_2$ \\
  & $\sfs{\discref}{(2,1)\ (11,4)}$ & 3 &
    $\Z \oplus \Z_2$ \\
  & $\sfs{\discref}{(2,1)\ (12,5)}$ & 3 &
    $\Z \oplus \Z_2 \oplus \Z_2$ \\
  & $\sfs{\discref}{(2,1)\ (13,5)}$ & 3 &
    $\Z \oplus \Z_2$ \\
  & $\sfs{\discref}{(3,1)\ (5,1)}$ & 3 &
    $\Z$ \\
  & $\sfs{\discref}{(3,1)\ (5,4)}$ & 3 &
    $\Z$ \\
  & $\sfs{\discref}{(3,1)\ (7,2)}$ & 3 &
    $\Z$ \\
  & $\sfs{\discref}{(3,1)\ (7,3)}$ & 3 &
    $\Z$ \\
  & $\sfs{\discref}{(3,1)\ (7,4)}$ & 3 &
    $\Z$ \\
  & $\sfs{\discref}{(3,1)\ (7,5)}$ & 3 &
    $\Z$ \\
  & $\sfs{\discref}{(3,1)\ (8,3)}$ & 3 &
    $\Z \oplus \Z_2$ \\
  & $\sfs{\discref}{(3,1)\ (8,5)}$ & 3 &
    $\Z \oplus \Z_2$ \\
  & $\sfs{\discref}{(4,1)\ (4,1)}$ & 3 &
    $\Z \oplus \Z_2 \oplus \Z_4$ \\
  & $\sfs{\discref}{(4,1)\ (4,3)}$ & 2 &
    $\Z \oplus \Z_2 \oplus \Z_4$ \\
  & $\sfs{\discref}{(4,1)\ (5,2)}$ & 3 &
    $\Z \oplus \Z_2$ \\
  & $\sfs{\discref}{(4,1)\ (5,3)}$ & 3 &
    $\Z \oplus \Z_2$ \\
  & $\sfs{\discref}{(5,2)\ (5,2)}$ & 3 &
    $\Z \oplus \Z_5$ \\
  & $\sfs{\discref}{(5,2)\ (5,3)}$ & 2 &
    $\Z \oplus \Z_5$ \\
  & $\sfs{\discref}{(2,1)\ (2,1)\ (3,1)}$ & 19 &
    $\Z \oplus \Z_2 \oplus \Z_2$ \\
  & $\sfs{\torus/o_2}{(3,1)}$ & 12 &
    $\Z \oplus \Z$ \\
  & $\sfs{\torus/o_2}{(3,2)}$ & 12 &
    $\Z \oplus \Z \oplus \Z_2$ \\
  & $\sfs{\kb}{(3,1)}$ & 21 &
    $\Z \oplus \Z$ \\
  & $\sfs{\kb}{(3,2)}$ & 21 &
    $\Z \oplus \Z \oplus \Z_2$ \\
  & $\sfs{\kb/n_3}{(3,1)}$ & 39 &
    $\Z \oplus \Z_{12}$ \\
  & $\sfs{\kb/n_3}{(3,2)}$ & 39 &
    $\Z \oplus \Z_2 \oplus \Z_6$ \\
  & $\sfs{\mobref}{(3,1)}$ & 15 &
    $\Z \oplus \Z$ \\
  & $\sfs{\mobref/n_2}{(3,1)}$ & 15 &
    $\Z \oplus \Z_{12}$ \\
  & $\sfs{\annrefref}{(3,1)}$ & 3 &
    $\Z \oplus \Z \oplus \Z_2$ \\
  & $\sfs{\annrefref/o_2}{(3,1)}$ & 3 &
    $\Z \oplus \Z_2 \oplus \Z_6$ \\
  & $\sfs{D}{(2,1) (2,1)} \bigcup_{\homtwotable{1}{1}{0}{1}} \sfs{\mb}{(2,1)}$ & 29 &
    $\Z \oplus \Z_8$ \\
  & $\sfs{D}{(2,1) (2,1)} \bigcup_{\homtwotable{-1}{2}{0}{1}} \sfs{\mb}{(2,1)}$ & 20 &
    $\Z \oplus \Z_2 \oplus \Z_2$ \\
  & $\sfs{D}{(2,1) (2,1)} \bigcup_{\homtwotable{0}{1}{1}{0}} \sfs{\mb}{(3,1)}$ & 10 &
    $\Z \oplus \Z_4$ \\
  & $\sfs{D}{(2,1) (2,1)} \bigcup_{\homtwotable{0}{1}{1}{0}} \sfs{\mb}{(3,2)}$ & 10 &
    $\Z \oplus \Z_2 \oplus \Z_2$ \\
  & $\sfs{D}{(2,1) (2,1)} \bigcup_{\homtwotable{1}{1}{0}{1}} \sfs{\annref}{(2,1)}$ & 8 &
    $\Z \oplus \Z_2 \oplus \Z_4$ \\
  & $\sfs{D}{(2,1) (2,1)} \bigcup_{\homtwotable{-1}{2}{0}{1}} \sfs{\annref}{(2,1)}$ & 6 &
    $\Z \oplus \Z_2 \oplus \Z_2 \oplus \Z_2$ \\
  & $\sfs{D}{(2,1) (2,1)} \bigcup_{\homtwotable{0}{1}{1}{0}} \sfs{\annref}{(3,1)}$ & 3 &
    $\Z \oplus \Z_4$ \\
  & $\sfs{D}{(2,1) (2,1)} \bigcup_{\homtwotable{0}{1}{1}{0}} \sfs{\annref}{(3,2)}$ & 3 &
    $\Z \oplus \Z_4$ \\
  & $\sfs{D}{(2,1) (3,1)} \bigcup_{\homtwotable{0}{1}{1}{0}} \sfs{\mb}{(2,1)}$ & 19 &
    $\Z$ \\
10 (ctd.)
  & $\sfs{D}{(2,1) (3,1)} \bigcup_{\homtwotable{-1}{1}{0}{1}} \sfs{\mb}{(2,1)}$ & 46 &
    $\Z$ \\
  & $\sfs{D}{(2,1) (3,1)} \bigcup_{\homtwotable{0}{1}{1}{0}} \sfs{\annref}{(2,1)}$ & 5 &
    $\Z \oplus \Z_2$ \\
  & $\sfs{D}{(2,1) (3,1)} \bigcup_{\homtwotable{-1}{1}{0}{1}} \sfs{\annref}{(2,1)}$ & 12 &
    $\Z \oplus \Z_2$ \\
  & $\sfs{D}{(2,1) (3,2)} \bigcup_{\homtwotable{0}{1}{1}{0}} \sfs{\mb}{(2,1)}$ & 19 &
    $\Z$ \\
  & $\sfs{D}{(2,1) (3,2)} \bigcup_{\homtwotable{0}{1}{1}{0}} \sfs{\annref}{(2,1)}$ & 5 &
    $\Z \oplus \Z_2$ \\
  & $\sfs{A}{(2,1)} / {\homtwotable{0}{-1}{1}{2}}$ & 41 &
    $\Z \oplus \Z_5$ \\
  & $\sfs{A}{(2,1)} / {\homtwotable{-1}{-2}{1}{1}}$ & 11 &
    $\Z \oplus \Z_6$ \\
  & $\sfs{A}{(2,1)} / {\homtwotable{1}{3}{0}{1}}$ & 4 &
    $\Z \oplus \Z_3$ \\
  & $\sfs{A}{(2,1)} / {\homtwotable{-1}{3}{0}{-1}}$ & 4 &
    $\Z \oplus \Z_3$ \\
  & $\sfs{A}{(3,1)} / {\homtwotable{0}{-1}{1}{0}}$ & 35 &
    $\Z$ \\
  & $\sfs{A}{(3,1)} / {\homtwotable{0}{-1}{1}{1}}$ & 7 &
    $\Z \oplus \Z_4$ \\
  & $\sfs{A}{(3,1)} / {\homtwotable{1}{2}{0}{1}}$ & 1 &
    $\Z \oplus \Z_2$ \\
  & $\sfs{A}{(3,2)} / {\homtwotable{0}{-1}{1}{0}}$ & 35 &
    $\Z \oplus \Z_2$ \\
  & $\sfs{A}{(3,2)} / {\homtwotable{0}{-1}{1}{1}}$ & 7 &
    $\Z \oplus \Z_5$ \\
  & $\sfs{A}{(3,2)} / {\homtwotable{1}{2}{0}{1}}$ & 1 &
    $\Z \oplus \Z_2 \oplus \Z_2$ \\
  & $\specialten$ & 13 &
    $\Z \oplus \Z_2$ \\
\end{supertabular}
\end{center}

\renewcommand{\arraystretch}{1}

%% file: paper.bbl
\newcommand{\noopsort}[1]{}
\providecommand{\bysame}{\leavevmode\hbox to3em{\hrulefill}\thinspace}
\providecommand{\MR}{\relax\ifhmode\unskip\space\fi MR }
\providecommand{\MRhref}[2]{%
  \href{http://www.ams.org/mathscinet-getitem?mr=#1}{#2}
}
\providecommand{\href}[2]{#2}
\begin{thebibliography}{10}

\bibitem{adams-gieseking}
Colin~C. Adams, \emph{The noncompact hyperbolic {$3$}-manifold of minimal
  volume}, Proc. Amer. Math. Soc. \textbf{100} (1987), no.~4, 601--606.

\bibitem{italian-nor6}
Gennaro Amendola and Bruno Martelli, \emph{Non-orientable 3-manifolds of small
  complexity}, Topology Appl. \textbf{133} (2003), no.~2, 157--178.

\bibitem{italian-nor7}
\bysame, \emph{Non-orientable 3-manifolds of complexity up to 7}, Topology
  Appl. \textbf{150} (2005), no.~1-3, 179--195.

\bibitem{regina}
Benjamin~A. Burton, \emph{Regina: Normal surface and 3-manifold topology
  software}, {\tt http://\allowbreak regina.\allowbreak sourceforge.\allowbreak
  net/}, 1999--2005.

\bibitem{burton-nor7}
\bysame, \emph{Structures of small closed non-orientable 3-manifold
  triangulations}, to appear in J. Knot Theory Ramifications, {\tt
  math.\allowbreak GT/\allowbreak 0311113}, November 2003.

\bibitem{burton-facegraphs}
\bysame, \emph{Face pairing graphs and 3-manifold enumeration}, J. Knot Theory
  Ramifications \textbf{13} (2004), no.~8, 1057--1101.

\bibitem{burton-regina}
\bysame, \emph{Introducing {R}egina, the 3-manifold topology software},
  Experiment. Math. \textbf{13} (2004), no.~3, 267--272.

\bibitem{burton-nor8}
\bysame, \emph{Observations from the 8-tetrahedron non-orientable census}, to
  appear in Experiment. Math., {\tt math.\allowbreak GT/\allowbreak 0509345},
  September 2005.

\bibitem{cormen-algorithms}
Thomas~H. Cormen, Charles~E. Leiserson, Ronald~L. Rivest, and Clifford Stein,
  \emph{Introduction to algorithms}, 2nd ed., MIT Press, Cambridge, MA, 2001.

\bibitem{cuspedcensusold}
Martin~V. Hildebrand and Jeffrey~R. Weeks, \emph{A computer generated census of
  cusped hyperbolic 3-manifolds}, Computers and Mathematics (Cambridge, MA,
  1989), Springer, New York, 1989, pp.~53--59.

\bibitem{closedhypcensus}
Craig~D. Hodgson and Jeffrey~R. Weeks, \emph{Symmetries, isometries and length
  spectra of closed hyperbolic three-manifolds}, Experiment. Math. \textbf{3}
  (1994), no.~4, 261--274.

\bibitem{0-efficiency}
William Jaco and J.~Hyam Rubinstein, \emph{0-efficient triangulations of
  3-manifolds}, J. Differential Geom. \textbf{65} (2003), no.~1, 61--168.

\bibitem{layeredlensspaces}
\bysame, \emph{Layered-triangulations of 3-manifolds}, Preprint, February 2006.

\bibitem{italian10}
Bruno Martelli, \emph{Complexity of 3-manifolds}, Preprint, {\tt
  math.\allowbreak GT/\allowbreak 0405250}, January 2005.

\bibitem{italian9}
Bruno Martelli and Carlo Petronio, \emph{Three-manifolds having complexity at
  most 9}, Experiment. Math. \textbf{10} (2001), no.~2, 207--236.

\bibitem{italian-decomp}
\bysame, \emph{A new decomposition theorem for 3-manifolds}, Illinois J. Math.
  \textbf{46} (2002), 755--780.

\bibitem{italianfamilies}
\bysame, \emph{Complexity of geometric three-manifolds}, Geom. Dedicata
  \textbf{108} (2004), no.~1, 15--69.

\bibitem{italianlists}
\bysame, \emph{Computer data on 3-manifolds}, {\tt http://\allowbreak
  www.\allowbreak dm.\allowbreak unipi.\allowbreak it/\allowbreak
  pages/\allowbreak petronio/\allowbreak public\_html/\allowbreak data.html},
  November 2004.

\bibitem{matveev-complexity}
Sergei~V. Matveev, \emph{Complexity theory of three-dimensional manifolds},
  Acta Appl. Math. \textbf{19} (1990), no.~2, 101--130.

\bibitem{matveev6}
\bysame, \emph{Tables of 3-manifolds up to complexity 6}, Max-Planck-Institut
  f\"{u}r Mathematik Preprint Series (1998), no.~67, available from {\tt
  http://www.\allowbreak mpim-bonn.\allowbreak mpg.\allowbreak de/\allowbreak
  html/\allowbreak pre\-prints/\allowbreak preprints.html}.

\bibitem{matveev11}
\bysame, \emph{Recognition and tabulation of three-dimensional manifolds},
  Dokl. Akad. Nauk \textbf{400} (2005), no.~1, 26--28.

\bibitem{matveev12}
\bysame, \emph{Tabulation of three-dimensional manifolds}, Russian Math.
  Surveys \textbf{60} (2005), no.~4, 673--698.

\bibitem{sedgewick}
Robert Sedgewick, \emph{Algorithms in {C}++}, Addison-Wesley, Reading, MA,
  1992.

\bibitem{snappea}
Jeffrey~R. Weeks, \emph{Snap{P}ea ({H}yperbolic 3-manifold software)}, {\tt
  http://www.\allowbreak northnet.\allowbreak org/\allowbreak weeks/\allowbreak
  index/\allowbreak SnapPea.html}, 1991--2000.

\end{thebibliography}
